\documentclass{article}

\usepackage{anysize}
\usepackage{amssymb,amsmath}
\marginsize{2cm}{2cm}{2cm}{2cm}
\usepackage[utf8]{inputenc}

\newcommand{\n}{\noindent}
\renewcommand{\epsilon}{\varepsilon}

\usepackage[dvipsnames]{xcolor}
\usepackage{amssymb,amsmath,mdframed,ytableau,tikz,extarrows,graphicx,float,enumerate,mathdots}
\usepackage[mathscr]{eucal}
\usepackage[nottoc]{tocbibind}
\usetikzlibrary{calc,arrows,math,snakes}
\usepackage[raggedright]{titlesec}
\usepackage[linktocpage=true,colorlinks=true]{hyperref}

\numberwithin{equation}{section}

\usepackage{amsthm}
\newtheorem{thm}{Theorem}[section]
\newtheorem{lem}[thm]{Lemma}
\newtheorem{cor}[thm]{Corollary}

\newtheorem{prop}[thm]{Proposition}
\newtheorem{defn}[thm]{Definition}

\newtheorem{rem}[thm]{Remark}

\title{The Poisson Geometry of Plancherel Formulas for Triangular Groups}
\author{Nicholas M. Ercolani\textsuperscript{$\dagger$}}
\date{}

\begin{document}

\maketitle

Department of Mathematics, University of Arizona, 617 N. Santa Rita Ave., Tucson, 85721-0089, AZ, USA

\begin{center}
    *Corresponding author. E-mail: ercolani@math.arizona.edu\\
    ORCID: 0000-0003-2010-4205\textsuperscript{$\dagger$}\\
\end{center}

\begin{abstract}
In this paper we establish the existence of canonical coordinates for generic co-adjoint orbits on triangular groups. These orbits correspond to a set of full Plancherel measure on the associated dual groups. This generalizes a well-known coordinatization of co-adjoint orbits of a minimal (non-generic) type originally discovered by Flaschka. The latter had strong connections to the classical Toda lattice and its associated Poisson geometry. Our results develop connections with the Full Kostant-Toda lattice and its Poisson geometry. This leads to novel insights relating the details of Plancherel theorems for Borel Lie groups to the invariant theory for Borels and their subgroups. We also discuss some implications for the quantum integrability of the Full Kostant Toda lattice.
\medskip

{\bf Keywords:} orbit method,  Dixmier-Pukanszky operator, Toda lattice, polarizations, invariant theory
\end{abstract}

\tableofcontents

\section{Introduction} \label{sec:intro}

The purpose of this paper is to amplify and further explore the deep connections between group representations and Poisson geometry (which is the geometry underlying Hamiltonian mechanics). The perspectives we bring to this are grounded in the theory of integrable systems especially as related to generalizations of the Toda lattice. More specifically we study the bridge between the Plancherel theorem for the real (upper) triangular matrix groups, i.e., those of the form (with $b_{ii}>0$)
\[
B =  \left[\begin{array}{ccccc}
b_{11} & * &\dots &\dots &*\\
& b_{22} & * & &\vdots\\
& & \ddots & \ddots &\vdots\\
& & & \ddots & *\\
& & & & b_{nn}
\end{array} \right]
\]
and the phase space geometry of generic coadjoint orbits for those groups, which is the phase space geometry underlying the full Kostant Toda lattice  \cite{bib:efs}. For groups and phase space geometries in general, this bridge is referred to as the {\em orbit method} program and so this work provides further explicit examples for that program, illuminating some recent developments there as well as showing new results.  
Though we restrict to the case of triangular groups, our presentation is framed to make it natural to extend our results to general Borel subgroups of $GL(n, \mathbb{R})$ as well as Borels for other real semisimple Lie groups and so we will often refer to the triangular group as a Borel subgroup. This work also opens a door for applications to quantum Toda lattices which will be discussed in our Conclusions. 
\smallskip

Before getting to the outline of our results we will provide some general motivation concerning Plancherel theory. (For more details on this background we refer the reader to \cite{bib:dm}.) Recall the original Plancherel theorem which is one of the cornerstones of classical harmonic analysis related to Fourier theory. In the case of the one-dimensional unitary group, $S^1$, Plancherel states that for a function $f(\theta) \in L^2(S^1)$ the map from $f$ to its Fourier coefficients, $\hat{f}(n) \in \ell^2(\mathbb{Z})$, is an isometry $\left(||f||_2 = ||\hat{f}||_2\right)$. Moreover, this isometry is explicitly given by 
\begin{eqnarray} \label{fourier}
f(\theta) = \sum_{n = - \infty}^\infty \hat{f}(n) e^{2\pi i n \theta}
\end{eqnarray}
in the sense that the series converges to $f(\theta)$ in the $L^2$-norm.  (In what follows we will refer to formulas like \eqref{fourier} as {\em Plancherel formulas}.) From the group theoretic perspective, \eqref{fourier} is seen as an expansion over all homomorphisms, $e^{2\pi i n \theta}$, indexed by $n \in \mathbb{Z}$, from $S^1$ to $\mathbb{C}^*$, which are usually referred to as {\em characters}. Since $S^1$ is a commutative group, one knows that  its irreducible representations are all 1-dimensional and these characters constitute all of the irreducible unitary representations of  $S^1$. The latter set, corresponding to the additive group $\mathbb{Z}$, is referred to as the dual group and denoted, in this case, by $\widehat{S^1}$. A similar story to the above holds for $L^2(\mathbb{R})$ where 
$\mathbb{R}$ is regarded as a group with respect to translation, corresponding to the Fourier integral representation  of $f(x) \in L^2(\mathbb{R})$. 

However the story changes significantly when considering harmonic analysis on a non-commutative group $G$. One dimensional homomorphisms no longer suffice to give a complete decomposition of $L^2(G)$. This was realized, even going back to Frobenius, in the case of finite groups, who understood that in the non-commutative setting there were irreducible representations (homomorphisms from $G$ to $GL(d, \mathbb{C})$) with $d > 1$, whose matrix coefficients needed to be involved in any Fourier-like decomposition of $L^2$. For compact groups the description of how this goes is based on the Peter-Weyl theorem \cite{bib:knapp}
which classifies all the irreducible unitary representations of $G$. Though no longer always one dimensional, these are still finite dimensional and constitute a discrete measure space, $\widehat{G}$,  of equivalence classes labelled by $\lambda$. For each $\lambda$, an element of the equivalence class is a matrix representation $\rho_\lambda(g)$ acting on a finite dimensional Hilbert space  $V_\lambda$. One still refers to 
$\widehat{G}$ as the dual group even though it will typically no longer be a group. In this setting the Fourier coefficient of $f(g) \in L^2(G)$ at representation 
$\lambda$ is a matrix gotten by integrating the adjoint of the representation matrix for $\lambda$ against $f(g)$ with respect to Haar measure on the group:
\[
\hat{f}(\lambda) = \int_G f(g) \rho_\lambda(g)^* d\mu_{\text{Haar}}(g)
\]
and the Plancherel formula is
\begin{eqnarray} \label{P-W}
f(g) = \sum_{\lambda \in \widehat{G}} \dim(V_\lambda) \langle \hat{f}(\lambda), \rho_\lambda(g)^* \rangle.
\end{eqnarray}
The contraction, $\langle \cdot, \cdot \rangle$ appearing here is defined by taking the trace of the product of the two matrices appearing in the pairing. The discrete measure, $\sum_{\lambda \in \widehat{G}} \dim(V_\lambda) \delta_\lambda$, is the measure on $\widehat{G}$, referred to as the Plancherel measure for $G$ and $\dim(V_\lambda) $ is the {\em multiplicity} with which the representation $\lambda$ appears as an isotypic (meaning equivalent) copy in the Plancherel formula.    
We also have the statement that this is an isometry,
\[
||f||^2_{L^2(G)} = \sum_{\lambda \in \widehat{G}} \dim(V_\lambda) \,\, || \hat{f}(\lambda)||^2.
\]
\medskip

The entry of Poisson geometry into this story may be seen as stemming from quantum mechanics.  We will again say just a few words of motivation here but refer the reader to \cite{bib:gust} for further details.  A quantum system with many symmetries typically has a Hamiltonian, given as a differential operator, that commutes with these symmetries. In cases where this leads to a maximal commuting set of observables (operators) one can expect to have a complete description of the coherent states of the system in terms of the spectra of these operators acting on their common eigenfunctions. This decomposition parallels that of the Plancherel formulae just discussed. For instance, in the case of $G = SO(3)$ these spectra are the quantum numbers of the hydrogen atom \cite{bib:singer}. The semiclassical limits of these quantum systems correspond to classical Hamiltonian dynamical systems in which the commuting operators correspond to Poisson commuting functions. In a series of papers from the '60s and '70s, Kirillov and Kostant \cite{bib:kir, bib:kostant3} introduced the idea of reversing this process thereby bypassing most of the the complications of quantum mechanics in favor of classical dynamical formulations on phase spaces built directly from the Lie algebra of the associated symmetry group of the system.  In a number of cases they were able to show that the orbits of the co-adjoint action of $G$ on the dual of its Lie algebra $\frak{g}^*$ are in one-to-one correspondence with the elements of 
$\widehat{G}$. (Physically these orbits are the level sets of energy Casimirs for the system.) In these descriptions the Liouville measure on $\frak{g}^*$ is related to the Plancherel measure on $\widehat{G}$. So it became of interest to extend this idea to more general groups which gave rise to the orbit method program mentioned at the outset. This has been an active area of research for many years now \cite{bib:kir2}. 

As a point of departure we make an aside here to note that for $S^1$ the Lie algebra is $i \mathbb{R}$. Its dual Lie algebra is comprised of the  linear functionals  $2\pi i \gamma \theta $. The induced action of conjugation on these functionals, which is  the co-adjoint action, is trivial since this group is commutative. So the co-adjoint orbits in this case are just the singleton linear functionals. However, due to the topological non-triviality of $S^1$, a condition is required to ensure that the characters, which are the exponentials of these functionals, be well-defined on $S^1$.  This requirement is that $\gamma = n \in \mathbb{Z}$. So, though this case is a bit degenerate, it is true that the elements of the dual group, $\widehat{S^1}$, are in one-to-one correspondence with appropriately {\em quantized} co-adjoint orbits. This quantization condition is the last vestige of the quantum mechanical connection in the orbit method. However, for the non-compact groups that we now turn to, no quantization conditions are in fact required. 

 \bigskip
 
In this paper we will be concerned with the case of non-compact, non-commutative groups, such as $B$. This presents new challenges including the fact that irreducible representations will now typically come in continuous series (analagous to Fourier integrals for $\mathbb{R}$) so that the Plancherel measure is not simply discrete and may not even be, prima facie, well-defined. Moreover, the irreducible representations will for the most part be infinite dimensional (and required to be unitary). Our particular interest is in Plancherel theorems for nilpotent and solvable groups. There has been progress here too, though only comparatively recently, thanks to the pioneering work of Moore, Wolf and their collaborators \cite{bib:mw, bib:lw}. The result of \cite{bib:lw} for the Plancherel formula of $B$ takes the form stated in Theorem \ref{B-Planch}. This result contains a description of the dual group $\widehat{B}$  in terms of parameters denoted 
$\frak{a}^*_\diamond$ with associated Plancherel measure as well as a new ingredient for which there is no classical analogue in \eqref{P-W}. This new ingredient is a {\em Dixmier-Pukanzsky operator}, denoted $D$. It is necessitated because, unlike the compact groups or even nilpotent groups, $B$ is not unimodular (i.e., its left and right Haar measures are not the same). We will review the derivation of this result and explain the meaning of all these terms in the first two sections of this paper. However, our goal is not simply to reproduce this result. Rather it is to understand the essential ingredients involved in a purely Poisson theoretic way that has clear significance for integrable systems theory. In particular we express $\frak{a}^*_\diamond$ explicitly in terms of the Casimirs that cut out the symplectic leaves of the Lie-Poisson structure associated to $B$. We also explicitly define the Dixmier-Pukanzsky operator for $B$ in terms of the invariant theory associated to the maximal unipotent subgroup of $B$. To our knowledge this is the first time this particular kind of link between representation theory and Poisson geoemetry has been made in the literature. These results then yield some novel insights into the integrable systems structure of more general Toda lattices.  This leads to the extension of a signature result from the classical tri-diagonal Toda lattice. We very briefly describe that now; the details are part of the primary content of this paper.
\smallskip

The Toda lattice \cite{bib:toda} is a dynamical system on $\mathbb{R}^{2n}$, with coordinates $(p_1,\ldots,p_n,q_1,\ldots,q_n)$. The system is Hamiltonian with respect to the standard symplectic structure on $\mathbb{R}^{2n}$ with Hamiltonian

\begin{equation} \label{hamiltonian}
H(p_1,\ldots,p_n,q_1,\ldots,q_n)=\dfrac{1}{2}\sum_{j=1}^n p_j^2 + \sum_{j=1}^{n-1}e^{q_j-q_{j+1}}.
\end{equation}
Under a remarkable coordinate transformation found by Flaschka  \cite{bib:fl}, \eqref{hamiltonian} can be
transformed to a Hamiltonian system with respect to the Lie-Poisson structure on the dual to the Lie algebra of $B$ (see Section \ref{sec:poiss} for a precise description of this Poisson structure). In natural coordinates on the dual algebra, the phase space for this system can be taken to be tridiagonal Hessenberg matrices; i.e., matrices $X$ of the form

\begin{equation} \label{hessenberg}
X = \left(\begin{array}{cccc} a_1 & 1\\ b_1 & a_2 & \ddots\\ &\ddots & \ddots & 1\\ &&b_{n-1}&a_n\end{array}\right).
\end{equation}
In these variables the Hamiltonian for the classical Toda lattice becomes
\[
H= \frac12 \text{Tr} X^2.
\]
This result is the signature result alluded to in the previous paragraph. There is a family of coadjoint orbits on this phase space parametrized by values of $\text{Tr} X$ which is an invariant of the coadjoint action (equivalntly a Casimir
of the Toda system). Corresponding to these orbits are infinite dimensional irreducible unitary representations of $B$ related to the action of $B$ on classical 
{\em Whittaker vectors}. This correspondence was worked out in detail by Kostant \cite{bib:kostant} as a significant illustration of the power of the orbit method in providing explicit descriptions of representations of solvable Lie groups. However, these orbits are among the smallest dimensional coadjoint orbits.
In this paper we initiate a similar analysis for the {\em largest dimensional} orbits.   These orbits correspond to a set of full measure in the dual group of $B$ which connects this analysis to a main theme of this paper concerning the Plancherel formula. Our approach may be understood as a kind of reversal of the approach taken by Flaschka. We will seek a coordinate transformation under which the full Kostant Toda lattice and its underlying Poisson geometry is expressed in terms of canonical coordinates with respect to a standard symplectic structure. We will prove the existence of such a transformation for the maximal dimensional (generic) 
orbits and some implications this has for the representation theory. We briefly describe some of those implications in the next paragraph.

\bigskip

The outline for this paper is as follows. In Section \ref{sec:back} we develop the necessary background on the construction of unitary irreducible representations for unipotent triangular matrices, in terms of a hierarchy of Heisenberg algebras. In Section
\ref{sec:borel} we develop the form of the Plancherel formula for $B$ that anchors our work. This entails developing some structures that make the Plancherel measure well-defined and get around a central problem related to the non-unimodularity of $B$. In particular this leads to the introduction of a Dixmier-Pukanszky operator. These two sections present material that already exists in the literature but which we summarize in some detail here in order to make our exposition self-contained.
The approach we take there follows more recent treatments given in \cite{bib:wo, bib:wo1}. In Section \ref{sec:poiss} we present the essential results on the Poisson structure of the dual Lie algebra for $B$ from \cite{bib:efs}.  In Proposition \ref{casimirs} the results on semi-invariants are presented differently and in more detail than in  \cite{bib:efs} as suits our needs here. We make use of more recent developments in invariant theory due to Kostant that appeared after publication of our previous work. We then derive, in Theorem \ref{N-invar} and Corollary \ref{DP-symbol}, the completely Poisson-theoretic  characterization of the Dixmier-Pukanszky operator associated with the Plancherel formula for $B$.  Section \ref{sec:polar} contains the remainder of our main results. In Theorem \ref{thm:puk} we determine a particularly natural polarization for the symplectic structure on the orbits enabling us to interpret these orbits as cotangent bundles on a $B$-homogeneous space, $T^*(B/H)$ where $H$ is the subgroup associated to the polarization.  The symplectic structure on the orbit then gets identified with the 
canonical symplectic structure of this cotangent bundle by a theorem of Pukanszky. In Theorem \ref{thm:puk-phase} the polarization is used to establish the existence of canonical coordinates, mentioned above, that interpolate between the symplectic leaf structure described in Section \ref{sec:poiss} on the one hand, and the Heisenberg hierarchy outlined in Section \ref{sec:back} on the other. Then in Section \ref{sec:Planch-rev} we relate all these Poisson perspectives back to the Plancherel formula for $B$ stated in Section \ref{sec:borel} and describe their significance. We will see that the Hilbert space for the representations corresponding to the generic orbits is naturally identified with $L^2(B/H)$. Finally, in Section \ref{conclusions} we outline some potential directions for future investigation that our work here may lead to. One of those directions is to return to the quantum systems with symmetries that gave rise to the orbit method in the first place. It will be of interest to consider quantum Toda lattices associated to the generic orbits. The results of this paper are a first step towards assessing whether or not these systems are quantum integrable.

\section{Background} \label{sec:back}
\subsection{Decomposition of a Real Semisimple Lie Groups} \label{sec:semidecomp}

In this section we will introduce the basic Lie group and algebra structures that we will be dealing with. 

 \n We consider the Lie algebra decomposition of $n \times n$ matrices
 
 \begin{eqnarray*} \label{frak}
 \frak{gl} = \frak{gl}(n, \mathbb{R}) & = & \frak{n}_- \oplus \frak{b}_+  
 \end{eqnarray*}
 where $\frak{n}_-$ is the lower triangular nilpotent sub-algebra (nilradical) and $\frak{b}_+$ is the complementary maximal solvable sub-algebra, referred to as a {\it Borel sub-algebra}:
 \begin{eqnarray*}
 \frak{n}_- =  \left(\begin{array}{ccccc}
0 &  & & &\\
*& 0 &  & &\\
\vdots & \ddots & \ddots & \ddots &\\
\vdots &  & \ddots & \ddots & \\
* & \dots & \dots & * & 0
\end{array} \right), &&
\frak{b}_+ = \left(\begin{array}{ccccc}
* & * &\dots &\dots &*\\
& * & * & &\vdots\\
& & \ddots & \ddots &\vdots\\
& & & \ddots & *\\
& & & & *
\end{array} \right).
\end{eqnarray*}
We will also use $\frak{n}$ to denote the transpose of $\frak{n}_-$ and $\frak{b}_-$ to denote the transpose of $\frak{b}_+$. Employing the {\it principal nilpotent} element, 
 \begin{eqnarray} \label{eps}
\epsilon &=& \left(\begin{array}{ccccc}
0 & 1 & & &\\
& 0 & 1 & &\\
& & \ddots & \ddots &\\
& & & \ddots & 1\\
& & & & 0
\end{array} \right)
\end{eqnarray}
we introduce the affine translate, 
\begin{eqnarray*}
\epsilon + \frak{b}_- &=& \left(\begin{array}{ccccc}
* & 1 & & &\\
*& * & 1 & &\\
\vdots & \ddots & \ddots & \ddots &\\
\vdots &  & \ddots & \ddots & 1\\
* & \dots & \dots & * & *
\end{array} \right), \\
\end{eqnarray*}
to represent the dual algebra $ \frak{b}^*_+ $ with respect to the $G$-invariant, non-degenerate inner product (Killing form), $(X,Y) = \text{Tr} XY$, on  $\frak{gl}$. 
The affine space $\epsilon + \frak{b}_- $ is the space of all lower Hessenberg matrices which we'll denote by $\mathcal{H}$.  It will be useful to introduce the following algebra and group projections,
\begin{eqnarray*}
 \pi_- : \frak{gl} \to  \frak{n}_-, \qquad
  \Pi_- : G \to  N_- \\
 \pi_+ : \frak{gl} \to  \frak{b}_+, \qquad
 \Pi_+ : G \to  B_+
\end{eqnarray*}
where $G = GL(n, \mathbb{R})$, $N_-$, the lower unipotent matrices, is the exponential group of the algebra $\frak{n}_-$ and $B_+$, the connected component of the identity in the invertible upper triangular matrices, is the exponential group of the algebra $\frak{b}_+$. $\Pi_\pm$ are defined on the open dense subset of $G$ where there is an $LU$ factorization. We also let $N$ denote the upper unipotent matrices.\\

\subsection{Decomposition of Unipotent Triangular subgroups} We introduce a further decomposition of the nilpotent algebra $\frak{n}$ which will play a fundamental role in the remainder of this paper.   Let $e_{i,j}$ denote the elementary matrix whose $(i,j)$-entry is 1 with all other entries being 0. Then we make the following definitions for $1 \leq r \leq \lfloor \frac{n}{2} \rfloor := R$.
\begin{eqnarray*}
\frak{m}_r &:=& \text{span} \{ e_{r,j} , e_{i,n-r+1} | r < i , j < n-r+1\} \cup \mathbb{R} e_{r,n-r+1}\\
\frak{n}_r &:=& \bigcup_{s \leq r} \frak{m}_s.
\end{eqnarray*} 
$\frak{m}_r$ is a sub-algebra of $\frak{n}$ isomorphic to the Heisenberg algebra of dimension $2(n - 2r) + 1$, for $0 < r \leq R$ and $\frak{n}_r $ is an ideal in 
$\frak{n}$. The center of $\frak{m}_r$ is one-dimensional:  $\frak{z}_r = \mathbb{R} \cdot e_{r,n-r}$. Note that when $n$ is even, $\frak{m}_R = \frak{z}_R$. 
We also define the following decompositions.
\begin{eqnarray} \label{lalg1}
\frak{m}_r &=& \frak{z}_r + \frak{v}_r  \\ \label{lalg2}
\frak{n} &=& \frak{s} + \frak{v} \,\, \text{where}\\ \label{lalg3}
\frak{s} &=& \oplus  \frak{z}_r \,\, \text{and} \,\, \frak{v} = \oplus  \frak{v}_r \\ \label{lalg4}
d_r &=& \frac12 \dim (\frak{m}_r/\frak{z}_r ) (= n - 2r) \implies \frac12 \dim (\frak{n}/\frak{s} ) = d_1 + \cdots + d_R
\end{eqnarray}

At the group level one has, correspondingly,
\begin{eqnarray} \label{group1}
N_r &=& M_1 M_2 \cdots M_r  \,\, \text{is normal in} \,\, N \,\, \text{and}\\ \label{group2}
N_r &=& N_{r-1} \rtimes  M_r\\ \label{group3}
N &=& M_1 M_2 \cdots M_R\\ \label{group4}
S &=& Z_1Z_2 \cdots Z_R = Z_1 \times \cdots \times Z_R \,\, \text{where} \,\, Z_r = \,\, \text{the center of} \,\, M_r.
\end{eqnarray}

\subsection{Schr\"odinger Representations of Heisenberg Groups} \label{schrod}
The representation theory of Heisenberg groups such as $M_r$, corresponding to the Heisenberg algebra $\frak{m}_r$ is well-established and enshrined in 
the Stone-von Neumann theorem. We consider the irreducible unitary representation $\pi_{\lambda_r}, \lambda_r \in \frak{z}^*_r \equiv \mathbb{R}$, of $M_r$, known as a {\it Schr\"odinger representation}, acting on the Hilbert space $L^2(\mathbb{R}^{d_r})$ and defined by
\[
\left[\pi_{\lambda_r} (\vec{q}, \vec{p}, c)f\right] (\xi) = e^{i \lambda_r (\vec{p}\cdot \xi + \frac12 \vec{p}\cdot \vec{q} + c )} f(\xi + \vec{q})
\] 
where $\vec{p}$ is the coordinate vector associated to the basis elements $e_{r,r+1}, \dots e_{r, r + d_r}$, $\vec{q}$ is the coordinate vector associated to the basis elements
$e_{r+1, n-r}, \dots e_{r + d_r, n-r}$ and $c$ is the coordinate associated to $e_{r,n-r+1}$. (Note that although these elements were introduced as associatted to a basis for the algebra $\frak{m}_r$, they also serve as coordinates on $M_r$ since this is an exponential group.)
\begin{thm} \cite{bib:foll}
Excluding $\lambda_r = 0$, up to unitary equivalence, the $\pi_{\lambda_r}$ are all the infinite dimensional irreducible unitary representations of $M_r$. 
\end{thm} 
\n From this, standard results on Plancherel's theorem for the Heisenberg group apply \cite{bib:mw}. The Plancherel measure on the unitary dual 
$\widehat{M_r}$ is supported on the $\pi_{\lambda_r}$ with $\lambda_r \ne 0$ and given by
$2^{d_r} d_r! |\lambda_r| d\lambda_r $. For $f \in L^1(M_r)$, define  
$\hat{f}(\lambda_r) = \int_{M_r} f(\vec{p}, \vec{q},c) \pi_{\lambda_r}(\vec{p}, \vec{q},c) d\vec{p}d\vec{q} dc  $. The version of Plancherel's theorem here that corresponds to the isometry of the Fourier transform mentioned in Section \ref{sec:intro} is given by 
\begin{eqnarray*}
||f||^2_{L^2(M_r)} = 2^{d_r} d_r! \int_{\frak{z}^*} ||\hat{f}(\lambda_r) ||^2_{HS} |\lambda_r|^{d_r} d \lambda_r
\end{eqnarray*}
for $f \in L^1(M_r) \cap L^2(M_r)$ and $||\cdot||_{HS}$ denoting the Hilbert-Schmidt norm. The result corresponding to \eqref{P-W} is
\begin{eqnarray*}
f(x) = 2^{d_r} d_r! \int_{\frak{z}^*} \Theta_{\pi_{\lambda_r}}(R_x f) |\lambda_r|^{d_r} d \lambda_r
\end{eqnarray*}
where $x$ is the group element corresponding to $(\vec{p}, \vec{q},c)$, $R_x$ denotes right translation acting on functions $\left( (R_x f)(g) = f(gx)\right)$ and
 $\Theta_{\pi_{\lambda_r}}$ is the distribution character of $\pi_{\lambda_r}$ defined by
 \begin{eqnarray*}
 \Theta_{\pi_{\lambda_r}}(f) = \text{Tr}\,\, \pi_{\lambda_r} (f) 
 \end{eqnarray*}
for $f \in C_c^\infty(M_r)$ .

\subsection{Stepwise Square Integrable Representations} \label{sec:stepwise}
The structure summarized in \eqref{lalg1} - \eqref{group4} sets one up to carry out a recursive construction of the irreducible unitary representations (referred to as {\em unirreps}) of $N$ from the Schr\"odinger
representations discussed in the previous subsection.  The full details of this construction are carried out in the foundational paper of Moore and Wolf \cite{bib:mw}
and subsequent work nicely reviewed in \cite{bib:wo1}. Here we will provide some background and then present those results needed for our work in this paper. 
 For ease of comparison we, for the most part, adopt the notation used in \cite{bib:wo1}. The basic inductive link for the recursion is \eqref{group2} : 
 $N_r = N_{r-1} \rtimes  M_r$. In the following subdivisions of this subsection we will first review the mechanism of real polarizations by which one lifts an irreducible representation from  $ N_{r-1} $ to $N_r$, then state the consequence of Mackey's {\em little group method} by which an irreducible unitary representation of $N_r$  is realized by the Hilbert space tensor of the lifted representation with one of the Schr\"odinger representations of $M_r$ described in Section \ref{schrod}. It is noted that this a set of full measure in $\widehat{N}_r$. Finally we present the explicit statement of Plancherel's formula for $N$.
 
 \subsubsection{Real Polarizations}
Consider a Lie algebra $\frak{g}$ with connected Lie group $G$.  In the matrix group case (which will be the case for us) $G$ acts naturally on $\frak{g}$ by conjugation and we denote this {\em adjoint action} by $Ad_g$  for $g \in G$ ($Ad_g(x) \doteq  g x g^{-1} \doteq g \cdot x $ for $x \in \frak{g}$). This action induces a corresponding {\em co-adjoint} action on linear functionals $\ell \in \frak{g}^*$; viz., $Ad^*_g(\ell)(x) \doteq \ell(Ad_g(x) ) \doteq g \cdot \ell$. Let $G_\ell$ denote the isotropy subgroup of $G$ at  $\ell$ under the co-adjoint action. 

\begin{defn} \label{pukn}
Fix a linear functional $\ell$ on $\frak{g}$. One says that a subalgebra $\frak{k} \subset \frak{g}$ is a {\em real polarization} at $\ell$ if 
\begin{enumerate}
\item $\left\{ Ad_g \frak{k}  | g \in G_\ell\right\} = \frak{k} $
\item $\ell([\frak{k} , \frak{k} ]) = 0$
\item $2 \dim \frak{k}  =  \dim \frak{g} + \dim \frak{g}_\ell$, where $\frak{g}_\ell$ is the Lie algebra of $G_\ell$. 
\end{enumerate}
\end{defn}
\noindent In doing calculations with co-adjoint actions for Lie algebras, $\frak{g}$, that are subalgebras of $\frak{gl}(n, \mathbb{R})$, it is usually most convenient to again use the invariant inner product, $(X,Y) = \text{Tr}XY$ to represent  $\frak{g}^*$ as we did in Section \ref{sec:semidecomp}. 
As a first example consider $\frak{g} = \frak{m}_r $. Here we will represent $\frak{m}^*_r $ in terms of the {\em opposite} algebra $\frak{m}^\dagger_r $, the transpose within $\frak{gl}(n, \mathbb{R})$. We consider the case where $\ell$ is given by $\ell(X) = (e_{n-r+1, r}, X)$ for $X \in \frak{m}_r$. Then the representation of the co-adjoint action on $\ell$ is given by $\pi_{\frak{m}^\dagger_r} Ad_g^{-1}(e_{n-r+1, r})$ for $g \in M_r$ where $\pi_{\frak{m}^\dagger_r}$ denotes projection into 
$\pi_{\frak{m}^\dagger_r}$. A straightforward calculation shows that $(M_r)_\ell = \mathbb{I} + \frak{z}_r $. We claim that a polarization at $\ell$ is given by 
$\frak{k}  = \frak{z}_r \oplus \frak{v}^+_r$ where $ \frak{v}^+_r = \text{span} \{ e_{r,j}  | r <  j < n-r+1\} $. From Definition \ref{pukn}, condition (1) is obvious, The bracket, $[\frak{k} , \frak{k}]  = 0$ and so condition (2) is immediate as well. Finally (3) holds since $\dim \frak{k}  = d_r +1$ while 
$\dim \frak{m}_r  + \dim (\frak{m}_r)_\ell = 2 d_r + 1 + \dim \frak{z}_r =  2 d_r + 1 + 1 = 2 \dim \frak{k} $. This establishes the claim. We note that in terms of the Schr\"odinger representations of $M_r$ described in Section \ref{schrod}, $ \frak{v}^+_r $ is coordinatized by $\vec{p}$ corresponding to the classical momenta of a free particle. Moreover, the Hilbert space for $\pi_{\lambda_r}$ is $L^2(M_r/K) = L^2(\mathbb{R}^{d_r})$ where $K$ is the connected group corresponding to 
$\frak{k} $ and $\mathbb{R}^{d_r}$ is coordinatized by $\vec{q}$ corresponding to the classical position variables  of a free particle. 
 
The construction of unirreps in $\widehat{N}$ will be built on a ladder of analogous polarizations. We now describe the first rung of that ladder from which steps onto the higher rungs will be clear.  So take $\frak{g} = \frak{n}_2$ and take $r=1$ in the construction of the previous paragraph so that $\frak{n}_1 = \frak{m}_1$, 
$\lambda_1$ be the linear functional  corresponding to $e_{n,1}$ with $\frak{w} _1$ denoting the polarization we found for this case. Extending it by zero, 
$\lambda_1$ is a linear functional on $\frak{n}_2$ still represented by $e_{n,1}$, but now regarded as an element of $\frak{n}^\dagger_2 $. To emphasize the distinction, this extension will be denoted by $\lambda_1'$. One may again calculate its co-adjoint isotropy to find that  $(N_2)_{\lambda_1'} = 
\mathbb{I} + \frak{z}_1 \oplus \frak{m}_2$. Proceeding as before one may establish that 
$\frak{w} _1' = \frak{z}_1 \oplus \frak{m}_2 \oplus \frak{v}^+_1$ is a polarization at  $\lambda_1'$. 
From this one sees that $N_2/K_1' = N_1/K_1$ and consequently $L^2(N_2/K_1' ) = L^2(N_1/K_1)$ which is the representation space of $\pi_{\lambda_1}$. It follows that this representation extends to a unirep, $\pi'_{\lambda_1}$ in $\widehat{N}_2$ such that $d\pi'_{\lambda_1}$ vanishes on $\frak{z}_2$.

\subsubsection{Constructing the Essential Part of $\widehat{N}$}
The key step for going forward is Mackey's Little Group Method which applies in the settting of inducing representations from a normal subgroup. This applies to the setting of the previous subsection to yield, for $N_2$, the following result. 
\begin{prop} \cite{bib:mackey} \label{mackey}
Elements of $\widehat{N}_2$ whose restrictions to $N_1$ are multiples of $\pi_{\lambda_1}$ are Hilbert space tensor products 
$\pi_{\lambda_1}' \widehat{\otimes} \, \gamma$ where $\gamma \in \widehat{M_2} = \widehat{N_2/N_1} $.
 \end{prop}
From Section \ref{sec:semidecomp} we know that the representations $\pi_{\lambda_2}$, associated to $ \lambda_2 \in \frak{z}^*_2$ corresponding to  
a non-zero multiple of $  e_{n-r+1, r} $, constitute a set of full measure in  $ \widehat{M_2}$. It follows from Proposition \ref{mackey} that the representations 
$\pi_{\lambda_1}' \widehat{\otimes} \,\pi_{\lambda_2}$, which we will denote by $\pi_{\lambda_1 + \lambda_2}$, are inequivalent unirreps in $\widehat{N}_2$. 
It is evident that this construction can be inductively continued to arrive at a family of unirreps on $\widehat{N} $ determined by a sequence of linear functionals 
$0 \ne \lambda_r \in  \frak{z}^*_r, r = 1, \dots, R$ denoted by $\pi_\lambda$ where $\lambda = \lambda_1 + \dots + \lambda_R$. The full description of properties related to Plancherel's formulae for $N$ are then given by the following fundamental result of Moore and Wolf.
\begin{thm} \cite{bib:mw}
Denote $\deg(\pi_\lambda) = \deg(\pi_{\lambda_1}) \dots \deg(\pi_{\lambda_R})$ with $\deg(\pi_{\lambda_r}) = |\lambda_r|^{d_r}$. Then the coefficients 
\\ $f_{z,w}(x) = \langle z, \pi_\lambda(x) w\rangle$ of the irreducible unitary representation $\pi_\lambda$ on $N$ are in $L^2(N/S)$ and satisfy 
\[
||f_{z,w}||_{L^2(N/S)}^2 = \frac{||z||^2 ||w||^2}{\deg(\pi_\lambda)} = \frac{||z||^2 ||w||^2}{\prod |\lambda_r|^{d_r}}.
\]
\end{thm}
\noindent Due to the recursive structure of the construction leading to these representations and the previous theorem, the $\pi_\lambda$ are referred to as 
the {\em stepwise square integrable} representations of $N$ relative to the decompositions \eqref{group1} - \eqref{group4}.
\begin{prop} \cite{bib:wo1}
Plancherel measure on $\widehat{N}$ (to be specified in Section \ref{N-Planch}) is concentrated on the stepwise square integrable representations $\pi_\lambda$. 
\end{prop}

\subsubsection{The Plancherel Formula for $N$} \label{N-Planch}
To now fully state the Plancherel formula for $N$ one would like to have a more explicit description of the distribution character for $\pi_\lambda$ as a tempered distribution. We start at the level of the Heisenberg algebra $\frak{m}_r$ on which one may regard $\lambda_r$ as a linear functional whose kernel is the non-central elements of $\frak{m}_r$. Consider the coadjoint orbit, $Ad^*(M_r)\lambda_r $, of $\lambda_r$ in $\frak{m}^*_r$. The distribution character, originally expressed as a trace, may be re-expressed as an integral over this orbit. One uses Lebesgue measure, $d\nu_r$ on $(\frak{m}_r/\frak{z}_r)^*$ normalized so that the Fourier transform is an isometry from $L^2(\frak{m}_r/\frak{z}_r)$ to $L^2(\frak{m}_r/\frak{z}_r)^*$. The tangent space to the orbit is an affine translate of a hyperplane in $\frak{m}^*_r$ and so $d\nu_r$ may be translated to a measure $d\nu_{\lambda_r}$ on the orbit. The orbital integral expressing the distribution character may then be written as
\[
\Theta_{\pi_{\lambda_r}}(f) = c^{-1}_r |\lambda_r|^{-d_r} \int_{Ad^*(M_r)\lambda_r} \widehat{f_1}(\xi) d\nu_{\lambda_r}(\xi)
\] 
where $c_r = 2^{d_r} d_r!$, $f_1(\xi) = f(\exp(\xi))$ and $\widehat{f_1}$ is the Fourier transform with respect to the Lebesgue measure on the orbit. 
This all extends to $N$ by defining 
\begin{eqnarray} \label{alg1}
\frak{t}^* &=& \{ \lambda = \lambda_1 + \cdots + \lambda_R \in \frak{s}^*, 0  \ne \lambda_r  \in \frak{z}_r^* \,\,  \forall r \}\\ \label{alg2}
\rho(\lambda) &=& \lambda_1^{d_r} \lambda_2^{d_r} \dots \lambda_R^{d_R}\\ \label{alg3}
c &=& c_1 c_2 \dots c_R = 2^{d_1 + \cdots + d_R}d_1! d_2! \dots d_R! \\ \label{alg4}
d\nu_{\lambda} &=& d\nu_{\lambda_1} \times \cdots \times d\nu_{\lambda_R} \\ \label{alg5}
\mathcal{O}(\lambda) &=& Ad^*(N)\lambda = Ad^*(M_1)\lambda_1 \times \cdots \times Ad^*(M_R)\lambda_R.
\end{eqnarray}
It is straightforward to check that $\frak{t}^*$ is a cross-section for the generic (maximal dimensional) co-adjoint orbits of $N$.
One then has
\begin{thm} \cite{bib:wo1} \label{nilplanch}
The distribution character of  $\pi_\lambda \in \widehat{N}$ is
\begin{equation} \label{nilplanch1}
\Theta_{\pi_{\lambda}}(f) = \text{Tr}\,\, \pi_{\lambda} (f) = \frac1c \frac1{|\rho(\lambda)| }\int_{\mathcal{O}(\lambda)} \widehat{f_1}(\xi) d\nu_{\lambda}(\xi)
\end{equation}
and $N$ has the Plancherel formula
\begin{equation} \label{nilplanch2}
f(x) = c \int_{\frak{t}^* } \Theta_{\pi_{\lambda}}(R_x f) |\rho(\lambda)|  d\lambda.
\end{equation} 
\end{thm}
\subsection{Poisson Geometry of Triangular Nilpotent Algebras}
From the construction in the previous subsection it is natural to consider 
\[
b_{\lambda_r} : = (x,y) \mapsto \lambda_r([x,y]) 
\]
which are symplectic bilinear forms, nondegenerate on $\frak{m}_r/\frak{z}_r$, respectively. Each defines a Poisson bracket for functions on  $\frak{m}^*_r$ referred to as the Lie-Poisson bracket and defined by $\left\{ F , G \right\} = b_{\lambda_r}(\nabla F, \nabla G)$.  With respect to this bracket the variables $\vec{q}, \vec{p}$ appearing in the Schr\"odinger representations of Section \ref{schrod} are canonical, corresponding to position and momentum. This all extends naturally to $N$ as 
\begin{eqnarray} \label{b-symp}
b_{\lambda} : = b_{\lambda_1} \oplus \cdots \oplus b_{\lambda_R} , 
\end{eqnarray}
a symplectic bilinear form,  nondegenerate on $\frak{n}/\frak{s}$ determining a Poisson bracket for functions on $\frak{n}^*$. The co-adjoint orbits, $\mathcal{O}(\lambda)$, appearing in \eqref{alg5} and \eqref{nilplanch1} are the symplectic leaves  for this Poisson structure, parametrized by $\frak{t}^*$. (Functions on $\frak{s}^*$ are the Casimirs of this Poisson structure.) This perfectly realizes the promise of the Kirillov-Kostant orbit method. The equivalence classes of irreducible unitary representations on a set of full measure in $\widehat{N}$
are in 1:1 correspondence with the generic co-adjoint orbits of $N$ (equivalently the generic symplectic leaves of the associated Poisson bracket), as expressed in the Plancherel formula of Theorem \ref{nilplanch}.  It is also of interest to note that $\rho(\lambda)$ introduced in the previous subsection equals the Pfaffian of $b_\lambda$.
\section{Borel Subgroups} \label{sec:borel}
The restricted triangular group, $N$, considered in the previous section, is a connected unipotent subgroup of the Borel subgroup, $B_+$ of real upper triangular matrices in $GL(n, \mathbb{R})$. We will restrict attention to the connected component of the identity in $B_+$ which is comprised of elements of $B_+$ whose diagonal entries are all positive. For simplicity of notation we will just denote this connected component by $B$. 

An obstacle arises in trying to extend, to $B$, the inductive types of constructions used in deriving the Plancherel formula for $N$; namely, whereas $N$
is unimodular, $B$ is not. A locally compact group, $Y$, is unimodular if its left invariant Haar measure equals its right invariant Haar measure; i.e., it has a 
bi-invariant measure. In the case of non-unimodular groups the two invariant measures are multiplicatively related by a function, $\delta_Y$ on the group, referred to as the modular function. Abelian groups are clearly unimodular. A group $Y$ is unimodular if and only if $Y/Z(Y)$ is, where $Z(Y)$ is the center of $Y$. Thus a straightforward induction on the length of the central series for the $n$-step nilpotent group $N$ establishes that $N$ is unimodular. However, $B$ is not and has a non-constant modular function that will be presented shortly. 

\subsection{Dixmier-Pukanszky Operator}

Non-unimodularity poses an issue for the formulation of a Plancherel formula.  Suppose relations such as those in Theorem \ref{nilplanch} were to hold for a non-unimodular group. Then for $x=1$, the identity element,  the corresponding version of \eqref{nilplanch2} would state that $f(1) = C \int_{\widehat{G}} \Theta_{\pi_{\lambda}}(f) d\mu_K$, where $\mu_K$ is either the left or right invariant measure. The left hand side is manifestly invariant under conjugation whereas, in general, conjugation transforms $\pi_\lambda(f)$, and hence \eqref{nilplanch1}, by the modular function. In the case of $B$ (and certain other parabolic subgroups) this mismatch can be compensated by the introduction of a semi-invariant operator, the Dixmier-Pukanszky operator, into the Plancherel formula. This is thanks to the work carried out by Wolf and collaborators \cite{bib:lw, bib:wo}, the result of which we'll now summarize. 

Note that $B=AN$, where $A$ is the split part of the maximal torus represented within in the diagonal matrices. The calculation of the modular function of $B$ is expedited by using the fact that 
\[
\delta_B(an) = \alpha(a) \delta_A(a) \delta_N(n)
\]
where $\alpha(a)$ is a positive function on $A$ such that for all $f \in C_c(N)$, $\int_N f(Ad_a(n)) d\mu_N = \alpha(a) \int_N f(n) d\mu_N$ and $d\mu_N$ is left Haar measure on $N$. Applying the change of variables formula one sees that $\alpha(a)$ is simply given by the inverse determinant of the Jacobian of the linear map 
$Ad_a|_N$. Since $\delta_A(a)= 1 = \delta_N(n)$,
\begin{eqnarray*} 
\delta_B(an) &=& \left( \det Ad_a|_N\right)^{-1}\\
&=& \prod_{i=1}^n a_i^{2i - n -1}
\end{eqnarray*}
for $a = \text{diag} (a_1, \dots, a_n)$. For later use we write this in a form more naturally related to the root structure of $\frak{gl}$,
\begin{eqnarray} \nonumber 
\delta_B(an) &=& \prod_{r=1}^R \left(\frac{a_r}{a_{n-r+1}}\right)^{d_r + 1}\\  \label{modB}
&=& \prod_{r=1}^R \left(\exp\beta_r{(\log a})\right)^{d_r + 1}
\end{eqnarray}
where $\{ \beta_1, \dots, \beta_R\}$ are the roots $\beta_r = e_{r,r}^* - e^*_{n-r+1, n-r+1}$.

A Dixmier-Pukanszky operator is a positive, self-adjoint invertible differential operator acting on $L^2(B)$, defined in terms of the Fourier transform and semi-invariant under the adjoint action of $B$. To define this for the triangular case Wolf \cite{bib:wo1} notes from  \eqref{alg1} and \eqref{alg2} that the Pfaffian 
$\rho(\lambda)$ may naturally be regarded as a polynomial on  $\frak{s}^*$, non-vanishing on $\frak{t}^*$. By duality one may also regard $\beta_r$ as a linear functional on $\lambda \in \frak{s}^*$ and then define $\text{Det}_{\frak{s}^*}(\lambda):= \prod_r \beta_r(\lambda) $. Then one can confirm the following by direct calculation. 

\begin{prop}  \label{DP}
 The product $\rho \cdot  \text{Det}_{\frak{s}^*}$ is an ${\rm Ad}(B)$-semi-invariant polynomial on $\frak{s}^*$ of degree $\frac12 (\dim \frak{n} + \dim \frak{s})$ and weight equal to the modular function $\delta_B$.
 \end{prop}
With $V = \exp(\frak{v})$ one defines the differential operator $D$ as the Fourier transform of this polynomial acting on $B = AN =AVS$ by acting on $S = \exp(\frak{s})$:
\begin{eqnarray} \label{DP1}
\mathcal{F}(D f)(\lambda, w) &=&  \rho(\lambda) \text{Det}_{\frak{s}^*}(\lambda) \mathcal{F}( f)(\lambda, w)\\ \label{DP2}
\mathcal{F}( f)(\lambda, w) &=& \int_S f(z, w ) e^{- i (\lambda, z)} dz
\end{eqnarray}
 where $z$ is a Euclidean coordinate on $S$ and $w$ is a coordinate on $AV$. Finally one has
 \begin{thm} \cite{bib:wo1} \label{DPop}
 $D$ is an invertible self-adjoint differential operator of order $\frac12 (\dim \frak{n} + \dim \frak{s})$ on $L^2(B)$ with dense domain of Schwarz class on $B$
 and is ${\rm Ad}(B)$ semi-invariant of weight equal to the modular function $\delta_B$. So $|D|$ is a Dixmier-Pukanszky operator on $B$.
 \end{thm}
 
 There is an alternative formulation of this $D$ that has a natural extension to more general minimal parabolic subalgebras developed in \cite{bib:lw} and \cite{bib:kostant4}. We will present this in Section \ref{sec:poiss} and reveal its interpretation in terms of the Casimirs for the Lie-Poisson structure on $\frak{b}_+$.

\subsection{Plancherel Theorem}

One can now formulate the Plancherel theorem for $B$.  Any $\lambda \in  \frak{t}^*$ has a non-zero projection on each summand  $\frak{z}_r^*$ of 
$ \frak{s}^*$. The action by conjugation  of $a \in A$ on $\frak{z}_r^*$ amounts to multiplication by $\exp(\beta_r(\log a))$. It follows that the isotropy 
algebra and subgroup of the coadjoint action of $A$ on $\frak{t}^*$ are, resepctively,
\begin{eqnarray} \label{diamond}
\frak{a}_\diamond &=& \{\xi \in \frak{a} | \beta_r(\xi) = 0, r = 1, \dots, R \}\\ \label{gp-diamond}
A_\diamond &=& \exp(\frak{a}_\diamond).
\end{eqnarray}
Note that these isotropies are the {\em same} at all $\lambda \in \frak{t}^*$. The dual group $\widehat{A_\diamond}$ is comprised of unitary characters, $\exp(i \phi(\log a))$
with $\phi \in \frak{a}^*_\diamond$.  Irreducible, unitary representations can then be constructed, as before, by Mackey's little group method. 
\begin{eqnarray} \label{IndRep}
\pi_{\lambda, \phi} = \text{Ind}_{NA_\diamond}^{NA} (\pi_\lambda' \otimes \exp(i\phi)); \,\, \lambda \in \frak{t}^*, \phi \in \frak{a}^*_\diamond
\end{eqnarray}
where $\pi_\lambda'$ is the natural extension of the stepwise square integrable representation $\pi_\lambda$ from $N$ to $NA_\diamond$.  
$\pi_{\lambda, \phi}$ and $\pi_{\lambda', \phi'}$ are equivalent if and only if $\lambda' \in {\rm Ad}^*(A) \lambda$ and $\phi = \phi'$. 
\begin{prop} \cite{bib:wo1} \label{supp-Planch}
Plancherel measure for $B = AN = NA$ is concentrated on the set of equivalence classes of the irreducible representations $\{ \pi_{\lambda, \phi}\}$ defined above. 
\end{prop}

\begin{thm} \cite{bib:wo1} \label{B-Planch}
The distribution character, $\Theta_{\pi_{\lambda}, \phi}(f):  f \mapsto \text{Tr}\,\, \pi_{\lambda, \phi} (f) $ is a tempered distribution. If $f$ is of Schwarz class on $B$,
then
\[
f(x) = c \int_{\frak{a}_\diamond^*} \int_{\frak{t}^*}  \Theta_{\pi_{\lambda}, \phi}(D (R_x f) |\rho(\lambda)| d\lambda d\phi.
\]
\end{thm}

We note that this realization of Plancherel's theorem reveals a relation between the Dixmier-Pukanszky operator $D$ constructed in  \eqref{DP1} and the Plancherel measure $d \mu =  |\rho(\lambda)| d\lambda d\phi$. In general neither of these is unique by themselves; however, $D \otimes \mu$ is unique. The choices made here are natural in terms of the constructions used. The nature of $D$  as a semi-invariant also suggests a connection to the universal enveloping algebra of $B$ which further suggests a relation to the Lie-Poisson structure on $\frak{b}_+^*$. This encourages one to expect a correspondence between $\widehat{B}$ as a measure space and the symplectic leaf (co-adjoint orbit) foliation for that Poisson structure with its associated Liouville measure. 

\section{Poisson Structures, Orbits and Invariant Theory} \label{sec:poiss}
In this section we develop the Poissson geometry that relates to the Plancherel theorem for $B$ that was presented in the previous section. This section culminates with results that give a Poisson semi-invariant description of the Dixmier-Pukanszky and related operators. We will also point out how this description connects the perspectives  of the representation theory community with those of the integrable systems community.

 \subsection{Poisson Reduction and Induction}
 
 As mentioned in Section \ref{sec:semidecomp}, we represent the dual space of $\frak{b}_+$, $\frak{b}^*_+$, with respect to the Killing form on $\frak{gl}$, as 
$\epsilon + \frak{b}_-$ which is just the space of lower Hessenberg matrices, $\mathcal{H}$. From a Lie theoretic perspective, this space is the translate of $ \frak{b}_-$ by the sum of the  simple positive root vectors, $\epsilon$. In this representation of the dual space, the Lie-Poisson bracket may be concretely expressed, for functions $f, g$ on  $\frak{gl}$, as 
\begin{eqnarray} \label{KKPB}
\{\tilde{f}, \tilde{g} \}(X) &=& \left( X, \left[ \pi_+ \nabla f(X),  \pi_+ \nabla g(X)\right]\right)
\end{eqnarray}
where $\tilde{f} = f|_{\mathcal{H}}$ and $\nabla$ denotes the gradient with respect to the Killing form. 
Determining the symplectic geometry underlying this structure is intimately related to the existence of commuting families of functions (Hamiltonians) with respect to the bracket and Casimirs (functions that Poisson-commute with {\em all} other functions). A method due to Thimm \cite{bib:thimm} provides a method of constructing maximal commuting families of Hamiltonians (also known as {\em completely integrable Hamiltonian systems}), in our case, from Casimirs of subalgebras of $\frak{gl}$ that contain $\frak{b}_+$, referred to as {\em parabolic subalgebras}. The method proceeds by considering inclusion chains of such parabolic subalgebras and the corresponding chain of projections between their duals. 

\begin{eqnarray} \nonumber 
\frak{b}_+ &=& \frak{p}_{m} \subset \cdots \subset \frak{p}_{k} \subset \cdots \subset \frak{p}_{0} = \frak{gl} \\ \label{parabolic}
\frak{b}_+^* &=& \frak{p}_{m}^* \leftarrow \cdots \leftarrow \frak{p}_{k}^* \leftarrow \cdots \leftarrow \frak{p}_{0}^* = \frak{gl}^*
\end{eqnarray}
The projections in \eqref{parabolic} are {\it Poisson maps},  meaning that the pullback of the Lie-Poisson bracket of two functions on 
 $\frak{p}_{k-1}^*$  equals the Lie-Poisson bracket bracket of their pullbacks on $\frak{p}_{k}^*$ . We let $P_k$ denote the parabolic group (which contains $B_+$) corresponding to $ \frak{p}_{k}$.
 
One can show that functions on  $\frak{p}_{k-1}^*$, invariant under the coadjoint action of $P_{k-1}$, Poisson commute with {\it all} functions on $\frak{p}_{k-1}^*$; i.e., they are Casimirs. Since the projections in (\ref{parabolic}) are all Poisson maps; one can collectively pull back all these sequentially invariant functions to 
$\frak{gl}^*$ to form a family of involutive (commuting) functions there. Then restricting this family to $\frak{b}_+^*$ gives an involutive family of functions on 
$\mathcal{H}$. For the details on all of this we refer the reader to \cite{bib:efs}. In that paper a particular chain of parabolics is chosen as well as a class of invariant functions on each parabolic yielding Casimirs. It is then shown that the resulting involutive family on $\mathcal{H}$ is maximal. The chain of parabolics used are those which are symmetric with respect to the anti-diagonal of the matrices in $\frak{gl}: \frak{p}_{r} = \left[\frak{g} \backslash \frak{n}_r\right]^\dagger, r = 0, \dots \left[\frac{n-1}{2}\right]$, where for a subalgebra $\frak{h} \subset \frak{gl}, \frak{h}^\dagger$ denotes the subalgebra which is the transpose of $\frak{h}$. The corresponding parabolic subgroup, $P_r$, is the set of elements in $G$ whose entries below the diagonal in the first $r$ columns and to the left of the diagonal in the last $r$ rows are zero. The maximal involutive family on $\epsilon + \frak{b}_-$ just described is realized in terms of the following construction.
The essential ideas underlying this construction were due to \cite{bib:dlnt} who used it to prove the complete integrability of the full symmetric Toda lattice and were adapted, in $\cite{bib:efs}$, to establish this for the full Kostant Toda lattice. We reproduce the main details of the latter setting here since that will be fundamental for what we do. 
\begin{defn} \label{chopmatrix}
For $r = 0, \dots, R$, and $X \in  \frak{gl}$, denote by $(X - \eta \mathbb{I})_{(r)}$ the result of removing the first $r$ rows and last $r$ columns of $(X - \eta \mathbb{I})$. We refer to this as the $r$-chop of $X$ and denote the coefficients in $\eta$ of its determinant as follows
\begin{eqnarray} \label{polychop}
\det (X - \eta \mathbb{I})_{(r)} = E_{0,r} \eta^{n -2r} + \cdots + E_{n - 2r,r}. 
\end{eqnarray}
\end{defn}
To make the connection to invariants we require an intermediate notion. We will again use the Killing form to make an identification of the dual algebra 
$\frak{p}_{r}^*$ with a translate of $\frak{p}_{r -}$, the transpose, in $\frak{gl}$, of  $\frak{p}_{r}$; specifically $\frak{p}_{r}^* \simeq \tau_r + \frak{p}_{r -}$ where 
$\tau_r = e_{1,2} + \dots + e_{r, r+1} + e_{n-r, n-r+1} + \dots + e_{n-1,n}$. Put differently, $\tau_r $ is the sum over the root vectors of the positive simple roots whose negatives do {\em not} belong to $\frak{p}_{r}$. Note also that $\tau_R = \epsilon$. For $ X \in \tau_r + \frak{p}_{r -}$, the co-adjoint action of  $P_r$ is realized by
\begin{eqnarray} \label{coad}
Ad^*_p X = \tau_r + \pi_{\frak{p}_{r -}} p^{-1} \cdot (X - \tau_r), \quad X \in \frak{p}_{r -}, \,\,\, p \in P_r
\end{eqnarray}
where $g \cdot Y = g Y g^{-1}$.

\begin{defn}
We will say that a function $f$ on $ \tau_r + \frak{p}_{r -}$ is a co-adjoint semi-invariant for $P_r$ if
\[
f\left(Ad^*_p X\right) = \chi(p)f(X) \quad \forall p \in P_r
\]
for some character $\chi$ which we refer to as the weight of the semi-invariant. When $\chi$ is the trivial character we say that $f$ is a co-adjoint invariant. Alternatively we may refer to it as a parabolic invariant or parabolic Casimir. 
\end{defn}
\begin{prop} \label{casimirs}
The $E_{m,r}(X)$ are all {\em semi-invariants} of the co-adjoint action of $P_r$ on $\frak{p}_{r}^*$, with weight
\begin{equation} \label{weight}
\chi_r(p) = \frac{p_{1,1} \cdots p_{r,r}}{p_{n-r+1, n-r+1} \cdots p_{n,n}}.
\end{equation}
\begin{enumerate}
\item[i)] For $r=0, \{ E_{m,0}\}$ are invariants under conjugation by $GL(n, \mathbb{R})$ and generate the full Poisson-commutative algebra of Casimirs on $\frak{gl}^*$.
This is equivalent to the family generated by  $\{ \text{Tr}X^m\}$  for $X \in \frak{gl}$. Upon restriction to $\frak{b}^* (X \in  \epsilon + \frak{b}_- = \mathcal{H})$, these are no longer all Casimirs but they do constitute an involutive family. $\text{Tr}X^2$ is the Hamiltonian for the famous Toda Lattice system--more precisely the Full Kostant Toda System \cite{bib:efs}. 
\item[ii)] Since, for each fixed $r$, the weight is common, one has collections of $P_r$-Casimirs
\[
I_{m,r} = E_{m,r}/E_{0,r}, \quad r = 1, \dots, \left[\frac{n-1}{2}\right]
\]
that are rational functions of the coordinates on $\frak{p}_{r -}$. Define the subset of {\em generic} matrices in $\mathcal{H}$ to be
\[
\mathcal{H}_{gen} = \left\{ X \in \epsilon + \frak{b}_- | E_{0,r}(X) \ne 0 \quad r = 1, \dots, \left[\frac{n-1}{2}\right]\right\}
\]
This is a Zariski open subset of  $\mathcal{H}$.
Taken all together, 
$\{ \text{Tr}X^m, I_{m,r}; r = 1, \dots,\left[\frac{n-1}{2}\right]\}$ is a maximal involutive family for the Lie-Poisson structure on  $\mathcal{H}_{gen}$ (though the elements of this family are not all $B$-invariamts). It is also a complete family of commuting constants of motion for the Full Kostant Toda flow, demonstrating the complete integrability of that system on generic co-adjoint orbits. 
\item[iii)] The leading coefficients $\{ E_{0,r}\}$ are invariants for the co-adjoint action of $N$ represented on $\frak{n}_-$. Note that here we allow $r$ to run from 
1 to $R$ and $R$ may be larger than $\left[\frac{n-1}{2}\right]$.
 \item[iv)] The subfamily generated by $\{ \text{Tr}X, I_{1,r}\}$ is a complete family of Casimirs for the Lie-Poisson bracket on $\mathcal{H}$. As such, the generic level sets of these Casimirs are the maximal dimensional symplectic leaves of the bracket and also the principal co-adjoint orbits of $B$. 
\end{enumerate}
\end{prop}
\begin{proof}
To establish that the $E_{m,r}(X)$ are indeed semi-invariants for $P_r$ with weight \eqref{weight} we begin by rewriting \eqref{polychop} in terms of a contraction 
on the exterior space $V_{n-r} = \bigwedge^{n-r} \mathbb{R}^{n}$. This is a fundamental representation space for $GL(n, \mathbb{R})$ given by
$$
g (v_1 \wedge \dots \wedge v_{n-r}) = g v_1 \wedge \dots \wedge g v_{n-r.}
$$
Since $GL(n, \mathbb{R})$ is a matrix subgroup of $\frak{gl}$ this definition extends to all elements of $\frak{gl}$. By restriction this also applies to 
$P_r, \frak{p}_{r }$ and $\frak{p}_{r -}$. 
With respect to the standard basis $e_1, \dots, e_{n}$ of $\mathbb{C}^{n}$
one defines a Hermitian inner product on $V_n$ by 
$$
\langle e_{i_1} \wedge \dots \wedge e_{i_n}, e_{j_1} \wedge \dots \wedge e_{j_n}\rangle = \delta_{i_1, j_1} \cdots \delta_{i_n, j_n}.
$$
Set $v^{(n-r)} = e_1 \wedge \dots \wedge e_{n-r}$ and $v_{(n-r)} = e_{r+1} \wedge \dots \wedge e_{n}$. These are, respectively, the highest and lowest weight vectors, with respect to lexicographic order, of this fundamental representation \cite{bib:fuha}. It is now straightforward to see that 
\begin{eqnarray} \label{repchop}
\det (X - \eta \mathbb{I})_{(r)} &=& \langle \left(X - \eta \mathbb{I}\right) v^{(n-r)},  v_{(n-r)} \rangle.
\end{eqnarray}
Applying the co-adjoint action \eqref{coad} to this representation one has
\begin{eqnarray} \label{ad1}
\langle \left(Ad^*_p X - \eta \mathbb{I}\right) v^{(n-r)},  v_{(n-r)} \rangle &=&
\left\langle \left(\tau_r + \left(\pi_{\frak{p}_{r -}} p^{-1} \cdot (X - \tau_r)\right) - \eta \mathbb{I}\right) v^{(n-r)},  v_{(n-r)} \right\rangle\\ \label{ad2}
&=& \left\langle \left( \left(\pi_{\frak{p}_{r -}} p^{-1} \cdot (X - \tau_r)\right) - \eta \mathbb{I}\right) v^{(n-r)},  v_{(n-r)} \right\rangle\\ \label{ad3}
&=& \left\langle \left( \left( p^{-1} \cdot (X - \tau_r)\right) - \eta \mathbb{I}\right) v^{(n-r)},  v_{(n-r)} \right\rangle\\ \label{ad4}
&=& \left\langle  p^{-1} \cdot \left(X - \tau_r - \eta \mathbb{I}\right) v^{(n-r)},  v_{(n-r)} \right\rangle\\ \label{ad5}
&=& \left\langle  p^{-1} \left(X - \tau_r - \eta \mathbb{I}\right) p v^{(n-r)},  v_{(n-r)} \right\rangle\\ \label{ad6}
&=& \left\langle   \left(X - \tau_r - \eta \mathbb{I}\right) v^{(n-r)},  (p^{-1})^\dagger v_{(n-r)} \right\rangle\\ \label{ad7}
&=& p_{1,1} \cdots p_{r,r} \det C_r \left\langle   \left(X - \tau_r - \eta \mathbb{I}\right) v^{(n-r)},  (p^{-1})^\dagger v_{(n-r)} \right\rangle\\ \label{ad8}
&=& \frac{p_{1,1} \cdots p_{r,r} \det C_r}{p_{n-r+1, n-r+1} \cdots p_{n,n}\det C_r}\left\langle   \left(X - \tau_r - \eta \mathbb{I}\right) v^{(n-r)},   v_{(n-r)} \right\rangle\\
\label{ad9}
&=& \chi_r(p) \left\langle   \left(X  - \eta \mathbb{I}\right) v^{(n-r)},   v_{(n-r)} \right\rangle.
\end{eqnarray}
where $C_r$ denotes the Levi factor of $p$ (the $(n-r) \times (n-r)$ central core of $p$). In passing from \eqref{ad1} to \eqref{ad2} we simply observed that 
$\tau_r$ is an element of $\frak{n}_r$. All elements of $\frak{n}_r$ annihilate $e_1$ and so also $v^{(n-r)}$. The difference between \eqref{ad2} to \eqref{ad3} is again an element of $\frak{n}_r$ and so these two lines are also equal. In  \eqref{ad4} the conjugation is simply passed through the identity and in  \eqref{ad5} that conjugation is written out explicitly. In \eqref{ad6} $p^{-1}$ is moved from the "bra" to the "ket" by transpose. The highest and lowest weight vectors are respectively stabilized by $p$ and $p^{-1}$ with multipliers as shown in \eqref{ad7} - \eqref{ad8} and which may be read of from the block form of $P_r$. Finally, as at the outset, $\tau_r$ may be removed to get \eqref{ad9}. This shows that \eqref{polychop} is a semi-invariant with weight \eqref{weight} for all $\eta$. Hence this is also the case for all $E_{m,r}(X)$.

For $r=0$ \eqref{polychop} is just the characteristic polynomial for $X \in \frak{gl}$ whose coefficients are clearly invariant under conjugation and generate the same algebra as $\{ \text{Tr}X^m\}$. The rest of statement (i) follows from the Adler-Kostant-Symes Theorem \cite{bib:efs}. 

The statement in (ii) about $P_r$ Casimirs is now immediate. The involutivity follows from Thimm's construction outlined earlier and the Adler-Kostant-Symes theorem. The maximality follows from the complete integrability on generic co-adjoint orbits of $B$ for which we refer the reader again to \cite{bib:efs}.  

Statement (iii) is of importance for what comes in the remainder of this section, As a consequence of the form of the weight \eqref{weight}, the $E_{m,r}(X)$ depend only on the $A$-factor of $B= AN$ and so they are $N$-invariant on $ \tau_r + \frak{p}_{r -}$. In general they will not be invariant when restricted to $\frak{n}_-$. 
However, since the unrestricted $E_{0,r}(X)$ only depend on coordinates in $\frak{n}_-$ the $N$-invariance persists for those coefficients. 

In a similar way one may observe that $E_{1,r}(X)$ depends only on coordinates in $\frak{b}_-$ and so it continues to be a $B$-semi-invariant when restricted to 
$\mathcal{H} (= \epsilon + \frak{b}_-)$. The independence of the $\{ \text{Tr}X, I_{1,r}\}$ on $\mathcal{H}_{gen}$ holds as in (ii) and the completeness of this family of Casimirs asserted in Statement (iv) is a consequence of the polarization to be described in Theorem \ref{thm:puk}.
\end{proof}
A key feature of these systems is that the invariants, $I_{m,r}$, are rational functions for $r > 0$. The pole locus, where the $I_{m,r}$ are undefined, corresponds to the locus in $\mathcal{H}$ where the co-adjoint orbits are not principal and so the corresponding symplectic leaves are not generic and one would not expect them to contribute non-negligible support to the Plancherel measure on  $\widehat{B}$. On the other hand these singularities could potentially pose problems for an orbit-theoretic formulation of Plancherel and its applications. As a first step toward addressing such issues we now bring in the invariant theory that plays a central role in the Lie theory and, particularly, in the representation theory discussed in Section \ref{sec:borel}. 

\subsection{Invariant Theory of $N$}

We consider, for the moment, the general setting of a Lie algebra $\frak{g}$ and its associated (connected) Lie group $G$. There is a linear isomorphism
\[
\sigma : \mathbb{S}(\frak{g}) \to \mathbb{D}(G)
\]
between $\mathbb{S}(\frak{g})$, the symmetric algebra over $\frak{g}$, and $\mathbb{D}(G)$, the space of all left-invariant differential operators, $D$, on $G$.
There is also a canonical isomorphism between $\mathbb{S}(\frak{g})$ and $\mathbb{S}(\frak{g}^*)$, the latter being the space of polynomials on $\frak{g}$. Then 
$\sigma$ essentially amounts to the correspondence between the differential operator and its symbol. We refer to Appendix \ref{app} for a precise description and details. Let $\mathbb{S}(\frak{g}^*)^G$ denote the co-adjoint invariant polynomials which, under the isomorphism $\sigma$, corresponds to the bi-invariant differential operators on $G$. One of the main points in this article is that this isomorphism lies at the heart of the relation between the Plancherel measure and the Dixmier-Pukanszky operator in our setting. 
\smallskip

We will now invoke a theorem of Kostant which describes these spaces of invariants in the cases of interest to us. To be more self-contained we will present this result with some detail and then apply it to the semi-invariant structures presented in Proposition \ref{casimirs}. Kostant begins by letting $\frak{n}_-$ represent $\frak{n}^*$ with respect to the Killing form and considers the co-adjoint action of $N$ which is represented as 
\begin{eqnarray*}
N \times \frak{n}_- &\to& \frak{n}_-\\
(n, \ell) &\mapsto& \pi_{\frak{n}_-}(Ad_{n^{-1}} \ell).
\end{eqnarray*}
One then restricts to the principal (maximal dimensional) orbits which comprise a Zariski open subset, $X$, of $\frak{n}_-$ (so $\overline{X} = \frak{n}_-$). One may further observe that the dual of $\frak{t}^*$, which we denote by  $\frak{t}^\times$,  provides a cross-section of this orbit. So $X = \cup_{\lambda \in \frak{t}^\times} \mathcal{O}(\lambda) $. Moreover, all the principal orbits, $\mathcal{O}(\lambda)$, are isomorphic to $N/S$. All this is consistent with what was described in \eqref{alg1} - \eqref{alg5}. It is also straightforward to check that $\frak{t}^*$ is an orbit under the co-adjoint action of $A$ and that, for $a \in A, 
a \cdot \mathcal{O}(\lambda) = \mathcal{O}(a \cdot \lambda)$. The co-adjoint action of $B = AN$ preserves $\frak{n}_-$ with the subgroup $N$ acting along its orbits while $A$ acts transversely as described in the previous sentence. 
\begin{thm} \cite{bib:kostant4} \label{thm:fibr}
One has an isomorphism
\[
X \to N/S \times   \frak{t}^\times.
\]
$X$ is an affine variety; letting $\mathcal{A}(X)$ denote its affine ring, one the has the following isomorphism of $B$-modules,
\[
\mathcal{A}(X) \simeq \mathcal{A}(N/S) \otimes \mathcal{A}(\frak{t}^\times).
\]
It follows from this that $\mathcal{A}(X)^N \simeq \mathcal{A}(\frak{t}^\times)$.
\end{thm}
The last observation of this theorem implies that the restricted action of the split torus $A$ on the module of $N$-invariants $\mathcal{A}(X)^N$ decomposes
as an action by $A$-weights on $\mathcal{A}(\frak{t}^\times)$.  We will refer to the following weight filtration.
\begin{eqnarray*}
\text{$A$-weight lattice} &\Lambda& \subset \frak{a}^*\\
&\cup& \\
\text{root lattice} & \Lambda_{ad} & \\
&\cup& \\
\text{cascade sublattice} &\Lambda_{\mathcal{B}} & = \oplus_{\beta \in \mathcal{B}} \mathbb{Z} \beta, 
\end{eqnarray*}
where $\mathcal{B}$ denotes the collection of roots $\{ \beta_1, \dots, \beta_R\}$ defined in \eqref{modB}. It is referred to as a {\em cascade of strongly orthogonal positive roots} because it is an example of Kostant's cascade construction \cite{bib:jo, bib:kostant4}. As a consequence of the fact that these roots are mutually orthognal,
the direct sum decomposition in the last line above follows. Thus, $\Lambda_{\mathcal{B}}$
is a free module of rank $R (= \lfloor \frac{n}2 \rfloor)$, in which each weight appearing occurs with multiplicty 1. 
If $M$ is a $B$-module  let $\Lambda(M) \subset \Lambda$ be the set of $A$-weights occurrirng in the weight decomposition of $M$. 
Since $\mathcal{A}(X)$ is an affine algebra one has natural inclusions
\[
\mathbb{S}(\frak{n}) \subset  \mathcal{A}(X) \subset \mathbb{Q}(\frak{n})
\]
where we note that $\mathbb{S}(\frak{n}) = \mathcal{A}(\frak{n}_-)$ and $\mathbb{Q}(\frak{n})$ is the quotient field of $\mathbb{S}(\frak{n})$. This induces inclusions on the $N$-invariants as $A$-modules.
\[
\mathbb{S}(\frak{n})^N \subset  \mathcal{A}(X)^N \subset \mathbb{Q}(\frak{n})^N.
\]
But it then follows from Theorem \ref{thm:fibr} that
\[
\Lambda\left(\mathbb{S}(\frak{n})^N\right) \subset  \Lambda(\mathcal{A}(\frak{t}^\times)) \subset \Lambda\left(\mathbb{Q}(\frak{n})^N\right).
\]
One sees from an explicit calculation that $\Lambda(\mathcal{A}(\frak{t}^\times)) = \Lambda_\mathcal{B}$ and so
\[
\Lambda\left(\mathbb{S}(\frak{n})^N\right) \subset \Lambda_\mathcal{B} \subset \Lambda\left(\mathbb{Q}(\frak{n})^N\right).
\]
Observe further that since $\mathbb{S}(\frak{n})$ is isomorphic to a polynomial ring, it is a unique factorization domain and so any element  
$q \in \mathbb{Q}(\frak{n})$ is uniquely expressible, up to a scalar multiple, as $q = f/g$ where $f$ and $g$ are relatively prime elements of 
$\mathbb{S}(\frak{n})$. By this uniqueness, if $q$ were an $N$ invariant, the action of $N$ would also need to preserve $f$ and $g$ respectively 
up to scalar multiples; i.e., up to respective characters. However, since $N$ is unipotent any such character is trivial and so $f, g \in \mathbb{S}(\frak{n})^N$. 
It follows that $\mathbb{Q}(\frak{n})^N$ is the quotient field of $\mathbb{S}(\frak{n})^N$ and by the last chain of inclusions one has
\begin{thm} \cite{bib:kostant4}
Every $A$-weight in $\Lambda\left(\mathbb{S}(\frak{n})^N\right)$ occurs with multipiicity 1 in $\mathbb{S}(\frak{n})^N$. Moreover, 
$\Lambda\left(\mathbb{S}(\frak{n})^N\right) = \Lambda_\mathcal{B} \cap \Lambda_{\text{dom}}$ where $\Lambda_{\text{dom}}$ is the set of dominant weights; i.e., those weights $\psi$ such that $(\psi, \phi) \geq 0$ for all positive roots $\phi$.   Every weight $\gamma \in \mathbb{Q}(\frak{n})^N$ occurs with multiplicity 1 
in $\mathbb{Q}(\frak{n})^N$ and has the form $\gamma = \nu - \mu$ for $\mu, \nu \in \Lambda\left(\mathbb{S}(\frak{n})^N\right)$. Additionally, the prime factors of any element in $\mathbb{S}(\frak{n})^N$ are also in $\mathbb{S}(\frak{n})^N$. 
\end{thm}

\subsection{The Dixmier-Pukanszky Operator Revisited}

In our case the generators in $\Lambda_\mathcal{B} \cap \Lambda_{\text{dom}}$ are of the form $\nu + \nu^*$ where $\nu$ is the highest weight of one of the fundamental representations and $\nu^*$ is the lowest weight or, alternatively, it may be characterized as the highest weight of the contragredient representation.
$\nu + \nu^*$ is then a dominant weight for the representation $V_\nu \otimes V^*_\nu$. But now one sees that the coefficients $E_{m,r}(X)$ are matrix elements
of such a representation; precisely that repreaentation for which $V_\nu = \bigwedge^{n-r} \mathbb{C}^{n}$. From this perspective one sees that it is most natural to view the $r$-chop  $(X - \eta \mathbb{I})_{(r)}$ appearing in \eqref{chopmatrix}  as inducing a family of elements in $\text{End}\,\, V_\nu =  V_\nu \otimes V^*_\nu $. All of the $E_{m,r}(X)$ are elements of this space corresponding to the weight $\nu + \nu^* \in \Lambda_\mathcal{B} \cap \Lambda_{\text{dom}}$; however, as a consequence of Proposition \ref{casimirs}(iii), {\em only one} of these is the multiplicity 1 representative of that weight in $\mathbb{S}(\frak{n}^*)^N$, namely 
$E_{0,r}(X)$. This establishes the following theorem. 

\begin{thm} \label{N-invar}
The invariant theory underlying the center of the enveloping algebras for $\frak{n}$  is explicitly described in terms of the $r$-chops as follows.
$\mathbb{S}(\frak{n}^*)^N \simeq \mathbb{R}\left[E_{01}(X), \dots. E_{0,R}(X)\right]$ where $E_{0,r}(X)$ is explicitly defined in Proposition \ref{casimirs}
with $X \in \frak{n}_-$.  $E_{0,r}(X)$ is homogeneous of degree $r$. 
\end{thm}

We now observe, from its modular weight \eqref{modB}, that the Dixmier-Pukanszky operator, is $N$-invariant; moreover, by Proposition \ref{DP} it depends only on $\frak{n}^*$ and so in fact is a bi-invariant differential operator in  $\mathbb{D}(N)$. Thus, its symbol, under $\sigma$ corresponds to an element of 
$\mathbb{S}(\frak{n}^*)^N$. It follows from the theorem that this symbol has a unique expression in terms of the $E_{0,r}(X)$.
\begin{cor} \label{DP-symbol}
The symbol of the operator, $D$, appearing in the Dixmier-Pukanszky operator for the Plancherel Theorem (Theorem \ref{B-Planch}) has a unique expression in 
$\mathbb{S}(\frak{n}^*)^N$ given by 
\begin{eqnarray*}
\sigma(D) &=& \prod_{r=1}^R E^2_{0,r}(X) \quad \text{for n odd}\\
&=& E_{0,R}(X) \prod_{r=1}^{R-1} E^2_{0,r}(X) \quad \text{for n even}
\end{eqnarray*}
and so is homogeneous of degree $R(R+1)$ for $n$ odd and degree $R^2$ for $n$ even. In either case, with $R = \lfloor n/2 \rfloor$ one sees that this degree
is equal to $\frac12\left( \dim \frak{n} +  \dim \frak{s} \right)$, coinciding with the order of the operator $D$ stated in Theorem \ref{DPop}.
\end{cor}
\begin{proof}
By the Theorem one knows that the symbol has a unique representation of the form
\[
\sigma(D) =  \prod_{r=1}^R E^{\alpha_r}_{0,r}(X) \quad \text{for unique} \,\, \alpha_r \in \mathbb{N}.
\]
From \eqref{weight} we know that the weight of $E_{0,r}$ is $\sum_{\ell=1}^r \beta_\ell$. Hence
\begin{eqnarray*}
\text{weight} \left( \sigma(D) \right) &=& \sum_{r=1}^R \alpha_r \left(\sum_{\ell=1}^r \beta_\ell \right)\\
&=& \sum_{\ell =1}^R  \left(\sum_{r = \ell}^R  \alpha_r \right) \beta_\ell 
\end{eqnarray*}
where in the second line one has just interchanged the order of summation.  On the other hand, by definition this weight is determined by \eqref{modB}, 
\[
\text{weight} \left( \sigma(D) \right) =  \sum_{\ell =1}^R (d_\ell + 1) \beta_\ell .
\]
Comparing these two expressions yields the system of equations
\[
\sum_{r = \ell}^R  \alpha_r  = d_\ell + 1
\]
which can be rewritten using \eqref{eps}, with $n$ replaced by $R$, in Toeplitz form as
\begin{eqnarray*}
(\mathbb{I} + \epsilon + \epsilon^2 + \dots + \epsilon^R) (\alpha_1, \dots, \alpha_R)^\dagger &=& (d_1 +1 , \dots, d_R + 1)^\dagger \\
(\mathbb{I} - \epsilon)^{-1} (\alpha_1, \dots, \alpha_R)^\dagger &=& (d_1 +1 , \dots, d_R + 1)^\dagger \\
(\alpha_1, \dots, \alpha_R)^\dagger &=&  (\mathbb{I} - \epsilon) (d_1 +1 , \dots, d_R + 1)^\dagger \\
&=&   (2 , \dots, 2)^\dagger \quad \text{for odd}\,\, n  \\
&=&   (2 , \dots, 2, 1)^\dagger \quad \text{for even}\,\, n 
\end{eqnarray*}
which establishes the stated expressions for $ \sigma(D)$. The evaluation of the homogeneous degree then follows from the respective $E_{0,r}$ having degree $r$.
\end{proof}

\subsubsection{Some History of Two Perspectives}

As may be seen from what has been discussed up to this point, the semi-invariants $E_{0,r}$ have played a significant role, separately, for developments both in representation theory of  Borel subgroups and in natural integrable systems on generic co-adjoint orbits of these Lie groups. However, to our knowledge, 
Corollary \ref{DP-symbol} represents the first time a connection between these two roles has been explicitly identified.  To help highlight this we briefly mention some historical antecedents in the separate disciplines. On the side of integral systems theory, the first reference we are aware of is due to Arhangel'skii \cite{bib:arh} who essentially identifies $E_{0,r}$ and $E_{1,r}$ by a direct inspection of minors as semi-invariants of the $B$-coadjoint action and applies the translation method of Miscenko-Fomenko to generate an involutive family.  In \cite{bib:trof}, Trofimov extends this approach to determine the dimension of symplectic leaves for other simple Lie algebras. Deift, Li, Nanda and Tomei \cite{bib:dlnt} independently studied the Lie-Poisson structure for  symmetric matrices that is associated to a different splitting of $\frak{gl}$. They introduce the chop semi-invariants in this symmetric case. This provides a more conceptual framework for generating semi-invariants and involutive families which can be related to the Ritz values of $X$.  The goal set in \cite{bib:efs} was to place the conceptual framework offered by the $r$-chops in its broadest possible setting by formulating it Lie theoretically. This was carrried out for the triangular group, as described earlier in this section, (with extensions indicated for other classical Borel subalgebras) in terms of the invariant theory of more general parabolic subgroups. This, in turn, was used to construct a toric linearization of the chop-invariant commuting flows on the associated flag manifolds of the semisimple Lie algebra. Gekhtman and Shapiro \cite{bib:gs} later extended this approach to establish complete integrability of the generic Toda flows for simple Lie algebras, including the exceptional ones. This work makes fundamental use of Lie-theoretic structure related to maximal cascades of strongly orthogonal roots, mentioned earlier, as developed by  Joseph and Kostant.

In representation theory, the semi-invariants $E_{0,r}$ appear first in Dixmier's early formulation \cite{bib:dix} of a Plancherel theorem for $N$. Then in the seminal paper \cite{bib:lw}, Lipsman and Wolf develop the central role of Dixmier-Pukanszky operators in more general Plancherel theorems. They relate this to a construction of $D$ in terms of semi-invariants for "good" parabolics and, in fact, see these arising from matrix coefficients with highest weights associated to cascade roots analogous to what was described just before Theorem \ref{N-invar} relating to the chop semi-invariants. Kostant developed a systematic analysis  of invariants in the enveloping algebra of $\frak{n}$ as we outlined earlier and subsequently \cite{bib:kostant5} used this to fill a gap in the Lipsman-Wolf construction. This did not alter anything in the conclusions of \cite{bib:lw} but it did make clearer to us the connection of that construction 
to our Poisson theoretic calculations in Proposition \ref{casimirs} and was the initial motivation for us to prove Theorem \ref{N-invar} and Corollary \ref{DP-symbol}.

\section{Polarization on Orbits and Canonical Structures} \label{sec:polar}
The connection between the symplectic geometry on the generic co-adjoint orbits of $B$ and the corresponding unirreps is further illuminated by determining canonical coordinates on these orbits which, intrinsically, means that the orbits may be represented as $T^*B/H$, the cotangent bundle of a homogenous space for an appropriate subgroup $H$. The means for doing this is to first find a real polarization (Definition \ref{pukn}) with respect to the symplectic structure on each orbit.  To do this in a uniform way with respect to the principal orbit foliation requires an additional condition which was formulated by Pukanszky \cite{bib:puk}.  

\subsection{Pukanszky Polarization for Borels} \label{sec:pukpolar}

Again making our usual dual identification for $\frak{b}^*$, choose an element $f \in \mathcal{H = \epsilon + \frak{b}_- }$ on a generic orbit and let 
$B_f$ denote its isotropy subgroup under the co-adjoint action.

\begin{defn} \label{puk}
We say that a subalgebra $\frak{h} \subset \frak{b}$ is a {\em real polarization} at $f$ satisfying Pukanszky's condition if
\begin{enumerate}
\item $\left\{ Ad_g \frak{h} | g \in B_f\right\} = \frak{h}$
\item $(f,[\frak{h}, \frak{h}]) = 0$
\item $2 \dim \frak{h} =  \dim \frak{b} + \dim \frak{b}_f$
\item $f + \frak{h}^\perp \subseteq \mathcal{O}(f)$, where $\frak{h}^\perp = \left\{  \ell \in \frak{b}^* | \ell(h) = 0, \,\, \forall h \in \frak{h} \right\} $
%=  \left\{ X \in  \epsilon + \frak{b}_- | (X, \frak{h}) = 0  \right\}  . 
\end{enumerate}
(4) is Pukanszky's condition.
\end{defn}
With such a real polarization in place one can then derive a canonical description of the co-adjoint orbits. This is done in terms of the standard moment mapping on $T^*B$, equipped with its canonical symplectic structure, that is associated to $B$ acting symplectomorphically by right translation. (We refer to \cite{bib:gust} for the background.) In terms of the trivialization of the cotangent bundle with respect to left invariant vector fields this takes the simple form 
\begin{eqnarray} \label{JR}
J_R: T^*B = B \times \frak{b}^* & \to & \frak{b}^*\\ \nonumber
(b, \ell) & \to & \ell.
\end{eqnarray}
Assume that $B_f$ is connected (which is the case for us). Then define $H$ to be the connected Lie subgroup of $B$ whose Lie algebra is $\frak{h}$. 
$H$ is a closed subgroup of $B$. We now also consider the moment map associated to the left action of $H$ which, in the trivialization, takes the form
\begin{eqnarray} \label{JL}
J_L: B \times \frak{b}^* & \to & \frak{h}^*\\ \nonumber
(b, \ell) & \to & \left(Ad^*_{b^{-1}} \ell \right)|_{\frak{h}},
\end{eqnarray}
which denotes the restriction of the coadjoint translate of $\ell$ to $\frak{h}$.
One is then in a position to carry out symplectic reduction with  respect to $J_L$, at a generic value $f$, yielding a reduced symplectic manifold
\begin{eqnarray} \label{fiber}
\left(  J^{-1}_L (f)/H , \omega_f \right)
\end{eqnarray}
where $\omega_f$ is the reduced symplectic form. The moment map $J_R$ restricts naturally to the fiber \eqref{fiber} to give a reduced map $J^f_R$ to $\frak{b}^*$. The existence and general form of canonical structures on the generic coadjoint orbits of $B$ will follow from the next proposition. We refer to \cite{bib:bgr} for general background and details on this fundamental result.

\begin{prop}  \label{cannonfiber}
The following statements are equivalent.
\begin{enumerate}
\item $\frak{h}$ satisfies Pukanszky's condition;
\item The symplectic reduction \eqref{fiber} may be naturally identified with the standard canonical structure  $\left( T^*(B/H), \omega_{can} \right)$.  The reduced momentum map $J^f_R: \left( T^*(B/H), \omega_{can} \right) \to \epsilon + \frak{b}_-$ is onto the symplectic leaf (co-adjoint orbit) through $f$, 
$\mathcal{O}(f)$;
\item The symplectic action of $B$ on $T^*(B/H)$  is transitive;
\item $J^f_R: \left( T^*(B/H), \omega_{can} \right) \to \left( \mathcal{O}(f),\omega_{\mathcal{O}(f)} \right)$ is a symplectic diffeomorphism 
where $\omega_{\mathcal{O}(f)}$ is minus the orbit symplectic form
\[
\omega_{\mathcal{O}(f)}(\ell)(X,Y) = - (\ell, [X,Y]), \quad \ell \in \mathcal{O}(f), \,\,\, X,Y \in \frak{b}.
\]
\end{enumerate}
\end{prop}

\begin{rem}
In general the reduced symplectic form involves a {\em magnetic term}; however, since we deal here with exponential groups, such a term will not appear.
\end{rem}

To describe the polarization, $\frak{h}$, that we will choose, we first take advantage of a cross-section for the generic coadjoint orbits of $B$ due to Arhangel'skii \cite{bib:arh} and also used to great effect by Gekhtman and Shapiro as a motivation in \cite{bib:gs}. It takes the following form.

 \begin{eqnarray} \label{basept}
 f  &=& \left(\begin{array}{ccccccc}
\kappa_1 & 1 &&&&&\\
& \kappa_2& 1 & &&&\\
& & \ddots & \ddots &&\\
& &  & \ddots &\ddots&&\\
& & \iddots&  &\ddots&1& \\
&1 & & &&\kappa_2& 1\\
1& & & &&& \kappa_1
\end{array} \right),
\end{eqnarray}
where the diagonal entries are symmetric across the anti-diagonal, the entries on the super-diagonal are all 1, the entries on the anti-diagonal below the diagonal are all 1 and all other entries are 0. 

We note that the parameters $\kappa_r$ appearing in \eqref{basept} are expressible in terms of the Casimirs $\{ \text{Tr}X, I_{1,r}\}$ (see Theorem \ref{thm:puk-phase}).  We also note that it is clear from the definition of the $\beta_r$ following \eqref{modB} that the  $\kappa_r, r = 1, \dots, n-R$ constitute coordinates on $\frak{a}_\diamond$. Hence the generic orbits form a Zariski open affine subvariety, $\mathcal{X}$, of 
$\epsilon + \frak{b}_- $ which can be presented as a fibration over $\frak{a}_\diamond$,
\[
\mathcal{X} = \cup_{a \in \frak{a}_\diamond} \mathcal{O}_a.
\]
 
We define the polarization with respect to $f$ to be the subspace
\begin{eqnarray} \label{h-basis1}
\frak{h} &=& \text{span}\left[ \{ e_{1,1} + e_{n,n}, \dots, e_{R,R} + e_{n - R+1, n-R+1}\} \cup  \{ e_{i,j} | i < j \leq n-i+1 \}\right] \quad \text{for $n$ even}; \\ \label{h-basis2}
&=& \text{span}\left[ \{ e_{1,1} + e_{n,n}, \dots, e_{R,R} + e_{n - R+1, n-R+1}, e_{R+1,R+1}\} \cup  \{ e_{i,j} | i < j \leq n-i+1 \}\right] \,\,  \text{for $n$ odd}.
\end{eqnarray}
To see that this is indeed a subalgebra of $\frak{b}$ we first bring together some notions and notations from sections \ref{sec:back} and \ref{sec:borel} (see in particular
\eqref{lalg2}, \eqref{lalg3}, \eqref{b-symp}, \eqref{modB} and \eqref{diamond}). 
\begin{eqnarray*}
\frak{a}_\diamond &=& \{\xi \in \frak{a}\, |\,  \beta_r(\xi) = 0, r = 1, \dots, R \}\\ %= \text{span} \{ e_{1,1} + e_{n,n}, \dots, e_{R,R} + e_{n - R, n-R}\}
\frak{s} &=& \oplus  \frak{z}_r = \text{span}\{ e_{1,n},  \dots, e_{r,n-r}, \dots, e_{R, n-R}\}\\
\frak{v}^+ &=& \text{span} \{ e_{i,j} | i < j < n-i+1 \} 
\end{eqnarray*}
where $\frak{v}^+$ may also be described as a maximal isotropic subspace of the symplectic form $b_\lambda$ on $\frak{n}/\frak{s}$. Thus
\begin{eqnarray} \label{h-decomp}
 \frak{h} = \frak{a}_\diamond  \oplus \frak{s} \oplus \frak{v}^+ .
 \end{eqnarray}
 It follows from the structure relations 
 \begin{eqnarray} \label{structure}
 [e_{i,j} , e_{k,\ell}] = \delta_{j,k} e_{i,\ell} -  \delta_{i,\ell} e_{k,j} 
 \end{eqnarray}
 for $\frak{gl}$ that
 \begin{eqnarray} \label{hh-decomp}
 [\frak{h}, \frak{h}] \subset \frak{v}^+ 
 \end{eqnarray}
 confirming that $\frak{h}$ is a subalgebra of $\frak{b}$. Next we calculate the isotropy subgroup of $f$. 
 \begin{lem} \label{f-isotrop}
 \[
 B_f  = A_\diamond.
 \]
 \end{lem}
 \begin{proof}
 Noting that $B = AN = NA$ it will suffice to consider the actions by $A$ and $N$ separately.  It is straightforward to check that 
 \begin{eqnarray*}
 Ad^*_a f &=& \epsilon + \pi_{\frak{b}_{-}} a^{-1} \cdot (f - \epsilon) \quad \text{for} \,\, a \in A_\diamond,\\
 &=& \epsilon +  a^{-1} \cdot (f - \epsilon)\\
 &=& \epsilon + f - \epsilon = f.
 \end{eqnarray*}
 Hence, $A_\diamond \subset B_f  $. On the other hand the action of an $a$ from the complement of $A_\diamond $ will non-trivially scale the anti-diagonal entries of $f$. Finally, for $n \in N$, $Ad^*_n f $ preserves the diagonal entries of $f - \epsilon $ since conjugation just spreads that element into $\frak{b}$. On the other hand there will be an element of $N$ preserving the anti-diagonal entries of $f - \epsilon $ if and only if there is an element $\xi \in \frak{n}$ whose bracket with the anti-diagonal part of $f - \epsilon $ is in $\frak{n}$. However, by direct inspection one sees that this bracket acts on $\xi $ (within $\frak{gl}$) by interchanging the top $R$ rows of $\xi$ with its bottom $R$ rows and (with a sign change) the last $R$ columns of $\xi$ with its first $R$ columns. In this way one sees that every entry of $\xi$ gets moved to at least one separate single entry in $\frak{b}_-$. Hence the bracket will be 0 only if $\xi = 0$. 
 \end{proof}
 We can now check that $\frak{h}$ satisfies the conditions of Definition \ref{puk} and deduce some corresponding consequences from Proposition \ref{cannonfiber}. 
 \begin{thm} \label{thm:puk}
 $\frak{h}$ given by \eqref{h-decomp} is a real polarization at $f$ specified by \eqref{basept}. Furthermore, 
 \[
 \left\{ Ad^*_h f | h \in H\right\} = f + \frak{h}^\perp
 \]
 which implies that Pukanszky's condition is satisfied.
  \end{thm} 
  \begin{proof}
  By Lemma \ref{f-isotrop}, condition 1 of Definition \ref{puk} becomes  the requirement that $\left\{ Ad_g \frak{h} | g \in A_\diamond \right\} = \frak{h}$. However, since $A_\diamond$ contains the identity this is immediate. It is straightforward to see that $(f, \frak{v}^+) = 0$. By \eqref{hh-decomp} it follows that condition 2 is also satisfied. Next, we calculate dimensions,
\begin{eqnarray*} 
\dim \frak{b}_f &=& \dim \frak{a}_\diamond \\
2 \dim \frak{h} &=& 2 \dim \frak{a}_\diamond  + 2 \dim \frak{s} + 2 \dim \frak{v}^+ \\
&=& 2 \dim \frak{a}_\diamond  + 2 \dim \frak{s} +  \dim \frak{n}/  \frak{s} \\
&=& 2 \dim \frak{a}_\diamond  +  \dim \frak{s} +  \dim \frak{n}\\
2 \dim \frak{h} - \dim \frak{b}_f &=& \dim \frak{a}_\diamond  +  \dim \frak{s} +  \dim \frak{n}\\
&=& n-R + R +  \dim \frak{n}\\
&=& n + \frac{n(n-1)}{2} = \frac{n(n+1)}{2} \\
\dim \frak{b} &=& \frac{n(n+1)}{2} ,
\end{eqnarray*} 
showing that condition 3 is satisfied. Finally, we address the Pukanszky condition. Applying duality with respect to the Killing from on $\frak{gl}$ it is straightforward to see that $\frak{h}^\perp = \left[\left(\frak{a}^\perp_\diamond \cap \frak{a}\right)  \oplus \frak{v}^-\right]^\dagger$.  We write
\[
f = \epsilon + \kappa + Z
\]
where $\kappa = \text{diag}(\kappa_1, \kappa_2, \dots, \kappa_2, \kappa_1)$ is the diagonal component of $f$ and $Z$ is the component whose entries on the anti-diagonal below the diagonal are all 1 with all other entries 0: $Z = e_{n,1} + e_{n-1,2} + \dots + e_{n-R+1, R}$. Recall once more that the co-adjoint action is given by
\[
Ad^*_h f = \epsilon + \pi_{\frak{b}_{-}} h^{-1} \cdot (f - \epsilon).
\]
We may consider the action on the components separately. For $\kappa$, since the $h$ act as raising operators, it is straightforward to see that 
$\pi_{\frak{b}_{-}} h^{-1} \cdot \kappa = \kappa$.  To evaluate the form of $\pi_{\frak{b}_{-}} h^{-1} \cdot Z$ it suffices to consider its infinitessimal form, 
$\pi_{\frak{b}_{-}} [\frak{h},Z]$; since we are dealing with exponential groups here, the exponential of the range of the infinitessimal action maps diffeomorphically onto the image of the co-adjoint action. By linearity, this reduces one to evaluating the bracket on basis elements in \eqref{h-basis1} or \eqref{h-basis2} for which
one may appeal to the structure relations \eqref{structure}. We know from Lemma \ref{f-isotrop} that $ [\frak{a}_\diamond, Z] = 0$ which implies that the exponential of this action is the identity on $Z$. Thus it suffices to just consider
 $[\frak{s},Z]$ and  $[\frak{v}^+,Z]$. In the first instance we evaluate 
 \begin{eqnarray*}
 && [ e_{r, n-r+1}, e_{n,1} + e_{n-1,2} + \dots + e_{n-R+1, R}] \quad \text{for} \,\, 1 \leq r \leq R\\
 &=& e_{r,r} - e_{n-r+1, n-r+1} \quad 1 \leq r \leq R
 \end{eqnarray*}
which spans $\left(\frak{a}^\perp_\diamond \cap \frak{a}\right)$.  In the second instance we consider 
\begin{eqnarray*}
 && \pi_{\frak{b}_{-}}  [ e_{i,j}, e_{n,1} + e_{n-1,2} + \dots + e_{n-R+1, R}] \quad \text{for} \,\, i < j < n-i+1\\
 &=& \pi_{\frak{b}_{-}} \left(e_{i, n-j+1} - e_{n-i+1, j}\right) \quad i < j < n-i+1\\
 &=& - e_{n-i+1, j} \quad i < j < n-i+1.
 \end{eqnarray*}
which spans $(\frak{v}^-)^\dagger$. Assembling all these observations we see that
\begin{eqnarray*}
\left\{ Ad^*_h f | h \in H\right\} &=& \epsilon + \pi_{\frak{b}_{-}} \left( \kappa + Z \oplus \left(\frak{a}^\perp_\diamond \cap \frak{a}\right) \oplus (\frak{v}^-)^\dagger\right)\\
&=& \epsilon + \kappa + Z \oplus \left(\frak{a}^\perp_\diamond \cap \frak{a}\right) \oplus (\frak{v}^-)^\dagger\\
&=& f + \frak{h}^\perp.
\end{eqnarray*}
This establishes the Theorem. 

\end{proof}

\subsection{Canonical Representation for the Poisson Structure of $\frak{b}^*$} \label{sec:canrep}

The real polarization just described now lends itself to a global coordinate description of $\frak{b}^*$ as a canonical phase space consistent with the Heisenberg algebras underlying the construction of the stepwise square integrable representations presented in Section \ref{sec:stepwise}. We can now give another presentation of the fibration of
\[
\mathcal{X} = \{ X \in \epsilon + \frak{b}_- \big| E_{0,r}(X) \ne 0, \forall r\}
\]
as an explicit product. We define
\[
\mathcal{F} = \{ f \in \epsilon + \frak{b}_- \big|  \,\, \text{such that $f$ is of the form} \,\, \eqref{basept}\}.
\]
We know from Proposition \ref{casimirs} (iv) that $\mathcal{X}$ is foliated by maximal symplectic leaves, $B/B_f$ and by Lemma  \ref{f-isotrop} each of these is diffeomorphic to $B/A_\diamond$. 
\begin{thm} \label{thm:puk-phase}
$\mathcal{F}$ is a cross-section of the symplectic leaf foliation and thus the map $\mathcal{X} \to \mathcal{F}$ which sends an element of $\mathcal{X}$ to its unique representative in $\mathcal{F}$, yield the diffeomorphism
\begin{eqnarray} \label{foliate}
 \mathcal{X} \to   B/A_\diamond \times \mathcal{F}.
\end{eqnarray}
\begin{enumerate} 
 
 \item[a)] The reduced symplectic manifold \eqref{fiber} at $f$ is explicitly realized by noting that
 \begin{eqnarray*}
 J^{-1}_L (f) &=& \left\{  (b, \ell) \in B \times (\epsilon + \frak{b}_-) \big| Ad_{b^{-1}} \ell \in f + \frak{h}^\perp \right\}\\
 &=& \left\{  (b, \ell) \in B \times (\epsilon + \frak{b}_-) \big| Ad_{b^{-1}} \ell \in f + \left[\left(\frak{a}^\perp_\diamond \cap \frak{a}\right)  \oplus \frak{v}^-\right]^\dagger \right\}
 \end{eqnarray*}
 where $ \frak{v}^- = \text{span} \{ e_{i,j} |  n-j+1 < i < j  \} $. Since $f$ is a regular value, $H$ acts freely and properly on this fiber, the quotient $J^{-1}_L (f)/H$ is a smooth manifold and $\omega_f$ is the unique symplectic form on this quotient whose pull-back to $J^{-1}_L (f)$ equals the restriction of the canonical form on $T^*B$ to $J^{-1}_L (f)$.   
 
 \item[b)] The tangent space to $B/B_f (= B/A_\diamond)$ at $f$ splits as
 \[
 T_f(B/B_f) = \left[\left(\frak{a}^\perp_\diamond \cap \frak{a}\right) \oplus \frak{v}^-\right] \bigoplus \left[ \frak{s} \oplus \frak{v}^+\right].
 \]
 
 Under the exponential map the first factor maps locally diffeomorphically into $B/H$ while the second factor is diffeomorphic to $H/B_f$.
 
 \item[c)] 
 $ \mathcal{X}$ has a natural coordinatization of the form 

\begin{eqnarray} \label{coordinate}
 \left(\begin{array}{ccccccc}
a_1 + q_{0,1 }& 1 & &&&&\\
p_{1,1} & a_2 + q_{0,2}& 1 & &&&\\
p_{2,1}& p_{1,2} & a_3 + q_{0,3}& 1 &&\\
p_{3,1}& p_{2,2} & p_{1,3} & a_4 &1&&\\
p_{4,1}&p_{3,2} & p_{0,3}& q_{1,3} &a_3 - q_{0,3}&1& \\
p_{5,1}& p_{0,2}&q_{1,2} &q_{2,2} &q_{3,2}& a_2 - q_{0,2}& 1\\
p_{0,1}&q_{1,1} & q_{2,1} &q_{3,1}  &q_{4,1} &q_{5,1} & a_1 - q_{0,1}
\end{array} \right)
\end{eqnarray}
in which the $q_{i,j}$ may be taken as local coordinates on $B/H$ near $f$ and the $p_{i,j}$ may be taken as local coordinates on the cotangent fiber to $B/H$ at $f$. There are, however, global constraints between these coordinates coming from the Casimirs on each orbit. 

 \end{enumerate}
\end{thm}
\begin{proof}
The map from $\mathcal{X}$ to $\mathcal{F}$ in \eqref{foliate} is explicitly determined by the relations \cite{bib:r},
\[
\kappa_1 + \dots + \kappa_r = \frac{\text{Tr} X - I_{1,r}(X)}{2}
\]
which, given the fixed values of the Casimirs, $I_{1,r}(X)$ and $\text{Tr}(X)$, on the orbit, determine the diagonal entries $\vec{\kappa}$ of $f$ in \eqref{basept}.
These Casimirs are independent on $\mathcal{X}$ and there are $\left[\frac{n+1}{2}\right]$ of them. Note also that $\vec{\kappa}$ runs from $\kappa_1$ to 
$\kappa_{\left[\frac{n+1}{2}\right]}$ and that $\dim \frak{a}_\diamond = n -R = \left[\frac{n+1}{2}\right]$. It follows from Lemma \ref{f-isotrop} that $\frak{a}_\diamond$ parametrizes the generic symplectic leaves. This is consistent with $f$ providing a cross-section of the coadjoint action and also confirms the statement in Proposition \ref{casimirs} (iv) that the $\{ \text{Tr} X, I_{1,r}\}$ are a maximal independent set of Casimirs for  $\mathcal{X}$. Again by Lemma  
\ref{f-isotrop} each $B/B_f$ is diffeomorphic to $B/A_\diamond$, which also gives the projection onto the first factor in \eqref{foliate}.
\smallskip

In part (a) of the theorem, for the initial equality, we first note that in \eqref{JL}, 
\[
T^*B = B \times \frak{b}^*  = B \times (\epsilon + \frak{h}),
\]
so that the moment map may be re-expressed, using duality, as
\begin{eqnarray*}
\left(Ad^*_{b^{-1}} \ell \right)|_{\frak{h}} &=& \pi_{\frak{h}^\dagger} \left(Ad_{b^{-1}} (\ell - \epsilon) \right) + \epsilon, 
\end{eqnarray*}
where $\ell$ on the RHS denotes the element in $\epsilon + \frak{b}_-$ dual, with respect to the Killing form, to the  linear functional $\ell$ on the LHS. 
The projection $\pi_{\frak{h}^\dagger} $ is projection onto the transpose of $\frak{h}$ which is contained in $\frak{b}_-$. From our arguments in the proof of 
 Theorem \ref{thm:puk} one sees that the kernel of $\pi_{\frak{h}^\dagger}$ within $\frak{b}_-$ is precisely $\frak{h}^\perp$. Hence the inverse image of $f$ under the moment map $J_L$ is given by the RHS of the first equality of (a). The second equality just uses the explicit determination of $\frak{h}^\perp$ found previously in the proof of Theorem \ref{thm:puk}. The remaining statements in (a) just follow from standard facts about symplectic reduction. 
\smallskip

Part (b) summarizes the Lie algebraic calculations made in Section \ref{sec:pukpolar}. The coordinatization described in part (c) follows from applying the exponential map locally to (b). 
\end{proof}

\begin{rem}
It should be noted that the coordinatization in \eqref{coordinate} reflects the $n$-step Heisenberg structure of the nilradical $N$ contained in $B$.
\end{rem}
  
\subsection{Plancherel's Formula Revisited} \label{sec:Planch-rev}

We recall from Theorem \ref{B-Planch}, in the statement of Plancherel's formula for $B$, that for $f(x)$ a Schwarz class function on $B$ one has
  
\[
f(x) = c \int_{\frak{a}_\diamond^*} \int_{\frak{t}^*}  \Theta_{\pi_{\lambda}, \phi}(D (R_x f) |P(\lambda)| d\lambda d\phi.
\]
Also, from \eqref{IndRep} we know that 
$\pi_{\lambda, \phi}$ and $\pi_{\lambda', \phi'}$ are equivalent if and only if $\lambda' \in {\rm Ad}^*(A) \lambda$ and $\phi = \phi'$. Using the Killing form on
$\frak{gl}$ we may identify elements $\lambda \in \frak{t}^*$ with anti-diagonal elements below the diagonal, none of whose entries vanish. As noted in the proof of Lemma \ref{f-isotrop}, conjugation by elements $a$ in the complement of $A_\diamond $ will non-trivially scale the anti-diagonal entries of the element representing $\lambda$. In fact given any other  such $\lambda'$ there is an $a$ that will move $\lambda$ to $\lambda'$. It follows that for fixed $\phi \in \frak{a}_\diamond^*$ all $\pi_{\lambda, \phi}$ are equivalent and so represent just one unirep equivalence class in $\widehat{B}$. Their distribution characters appear in the Plancherel formula as isotypic summands in an integral direct sum, which one may regard as a representation of their multiplicity in the Plancherel measure. 
On the other hand the $\pi_{\lambda, \phi}$ are inequivalent for different $\phi$. It follows from  Proposition \ref{supp-Planch} that these parametrize a set of distinct equivalence classes comprising a set of full measure in $\widehat{B}$. 

We may also view this from the Poisson perspective. From Theorem \ref{thm:puk-phase} we know that the principal co-adjoint orbits in $\frak{b}^*$
are parametrized by the diagonal parameters $\vec{\kappa} = (\kappa_1, \dots, \kappa_{n - R})$ appearing in $f$ \eqref{basept}. However, under duality $\vec{\kappa}$ corresponds 1:1 with
$\frak{a}_\diamond^*$. So we have
\begin{thm}
The just described  correspondence between an open dense subset of $\widehat{B}$ and generic co-adjoint orbits in $\frak{b}^*$ explicitly verifies the validity of the orbit method that was abstractly established by Auslander and Kostant \cite{bib:ak}.  
Moreover, in direct relation to Proposition \ref{cannonfiber} (2) we may identify the Hilbert space for the representation $\pi_{\lambda, \phi}$ with $L^2(B/H)$. Also the symbol of the Dixmier-Pukanszky operator, $D$ appearing in the Plancherel formula is explicitly given by Corollary \ref{DP-symbol}. 
\end{thm}

\section{Conclusions} \label{conclusions}
In \cite{bib:ak} Auslander and Kostant established, on general grounds, that for type I solvable Lie groups there is a one-to-one correspondence between equivalence classes of its unireps, the dual group, and its co-adjoint orbits. Kostant \cite{bib:kostant} went on to make this explicit in the special case of the so-called Whittaker representations. These correspond to the {\em smallest} indecomposable co-adjoint orbits which are the tridiagonal (Jacobi) orbits in the case of triangular groups. That work opened the door for deep connections between representation theory and integrable systems theory related to the classical Toda lattices. More recently this has led to important relations with probability theory \cite{bib:bbo, bib:o}, number theory \cite{bib:sh} and combinatorics \cite{bib:bbf}. 

The results in this paper lay the groundwork for the kind of explicit description that Kostant carried out but for generic elements of the dual group, those corresponding to the {\em largest} dimensional co-adjoint orbits. Specifically we have made  explicit connections between all the elements of the Plancherel formulae
for triangular groups, $B$,  as developed in the works of \cite{bib:lw, bib:wo} on the one hand and the Poisson geometry of co-adjoint orbits on $\frak{b}^*$ as developed within the framework of the full Kostant Toda lattice. This has been done in a way that will make it natural to extend our results to the Borel subgroups of general complex semi-simple Lie groups.  This now positions one to explore deeper connections with integrable systems and allied fields similar to what has been seen for the Whittaker-classical Toda case. We expect this will be pursued in future investigations, but for now we will just take some space below to point out how the work in this paper may specifically contribute to those investigations.

\subsection{Canonical Coordinates and Standard Hamiltonian Systems} \label{sec:canco}

As with the construction of polarizations and canonical structures for generic orbits described in Section \ref{sec:polar}, one may carry out an analogous construction for the minimal orbits of tridiagonal form. We refer the reader to  \cite{bib:bgr} for a detailed description of this. Under the corresponding coordinate transformation the tridiagonal Toda system determined by \eqref{hessenberg} with Hamiltonian $TrX^2$ takes the form \eqref{hamiltonian} in terms of 
{\em standard} canonical variables $\{ p_1, \dots, p_n; q_1, \dots, q_n\}$.
By "standard"  here one means that in these variables the Poisson bracket ``diagonalizes" as $ \{p_i, p_j\} = 0 =  \{q_i, q_j\} , \{p_i, q_j\} = \delta_{i,j}$.

As a consequence of of Theorem \ref{thm:puk}, Proposition \ref{cannonfiber} implies that such coordinates exist in the case of generic orbits; however, it is not immediately clear what the explicit form of the coordinate transformation is or what form the Toda Hamiltonian, or any of the other Hamiltonians that commute with it, will take. However, in \cite{bib:symes} Symes outlines a general method for constructing the transformation based on the polarization. Moreover, the stratified Heisenberg structure evident in \eqref{coordinate} suggests that the canonical coordinates we seek may have a close relation to canonical Heisenberg variables. So this is a promising beginning that will be pursued in a future study.

\subsection{Quantization of Invariants and Quantum Integrability}

As mentioned in the introduction, the orbit method had its origins in the quantum mechanics of systems with many symmetries. Having worked out the details of this method in the case of triangular groups it is natural to go back and ask what implications this might have concerning associated quantum systems and in particular the quantum analogue of the Full Kostant Toda lattice. In the case of the tridiagonal orbits this has been done for the original Toda lattice. In this
case the quantum Toda lattice takes the form
\begin{equation} \label{quantham}
\widehat{H} = - \dfrac{\hbar^2}{2}\sum_{j=1}^n \frac{\partial^2}{\partial q_j^2} + \sum_{j=1}^{n-1}e^{q_j-q_{j+1}}
\end{equation}
following the formal prescription of canonical quantization which replaces $p_j$ by $i \hbar \partial/\partial q_j$ . 
A number of obstacles arise in trying to place this formal procedure on a rigorous footing. The first potential issue is that the choice of quantization is in general not unique; there may be several choices of operator whose semi-classical limit leads to the same $H$. (See \cite{bib:foll} for further details on this issue.) However, in the case of  \eqref{hamiltonian} the evident separation of variables obviates any such ambiguity. The relation between $H$ and \eqref{quantham} is consistent with the symbol map as described in Appendix \ref{app}.
%with $\widehat{H}$ regarded as a operator in $L^2{B/K}$ for the polarization group, $K$, associated to the tridiagonal orbit.
The next potential obstacle concerns the extension of classical integrability, and in particular Poisson commutativity, to the quantum level.  For the tridiagonal case we may consider the Poisson commuting invariants given by $\{ E_{m,0}(\vec{p}, \vec{q})\}$in Proposition \ref{casimirs} (i). These are the coefficients of 
$\det (X - \eta \mathbb{I})$ for $X$ of the form \eqref{hessenberg} when $a_i = p_i$ and $b_i = \exp{q_i - q_{i+1}}$. We note that when restricted to this tridiagonal phase space all the other semi-invariants, $\{ E_{m,r}\}$ for $r > 0$, that were mentioned in Proposiiton \ref{casimirs} vanish identically. However, the invariants for $r = 0$ suffice for the complete integrability of tridiagonal Toda. Once again, remarkably, the quantization of these invariants, denoted $\widehat{E}_{m,0}$ are uniquely determined, gotten again by simply replacing the  $p_j$ by $i \hbar \partial/\partial q_j$ \cite{bib:gi}. However, Poisson commutativity does not necessarily imply that these associated operators commute since multiplication of symbols does not pass, mutatis mutandis, over to composition of operators. The only thing one may conclude, without further examination, is that commutators vanish to order $\hbar^2$ \cite{bib:foll}. Nevertheless one {\em does} find that
$[\widehat{H}, \widehat{E}_{m,0}] = 0$ for all $m$, meaning that these operators play the role of quantum conservation laws. They are also referred to as a {\em complete set of commuting observables}, and they all have common eigenstates. This is what is meant by quantum integrability. The role of the common eigenstates is what is played by the Whittaker functions alluded to at the start of these Conclusions. 

In trying to extend quantization and related matters to the setting of the full Kostant Toda lattice all the above mentioned obstacles rear their heads again. Since the Hamiltonian and related invariants will be evaluated on generic orbits rather than just tridiagonal ones (and will now involve all the $\{ E_{m,r}\}$) what made things work out in the tridiagonal case may no longer be available. The passage to canonical variables based on our results here and envisioned in Section \ref{sec:canco} may be of decisive help, but that remains to be determined. There is one further potential obstacle in this extension to quantum full Toda, for which the results of this paper may help. As we saw, the invariants for $r>0$ are necessarily rational functions whose associated quantized operators should be given by pseudo-differential operators in terms of an appropriate Fourier transform. To formulate this we can make use of the Plancherel formula for
$B$ (Theorem \ref{B-Planch}) and define, using the invariance of  $I_{m,r}$,
\[
\widehat{I}_{m,r}f(x) = c \int_{\frak{a}_\diamond^*} \int_{\frak{t}^*}  \Theta_{\pi_{\lambda}, \phi}(D \widehat{I}_{m,r}  (R_x f) |\rho(\lambda)| d\lambda d\phi
\]
where $D \widehat{I}_{m,r}$ is a {\em rational} semi-invariant operator with the same modular weight as $D$ and defined in the same manner as $D$ was in 
\eqref{DP1} - \eqref{DP2} with symbol $\sigma\left( D \widehat{I}_{m,r} \right)$. (We note that the linear bijection  \eqref{symbol} extends to the function fields of the enveloping algebras involved.) The concern here is that the pseudo-differential representations of the quantizations of $I_{m,r}$ might involve singular integrals. By the results of this paper, and in particular Corollary \ref{DP-symbol}, one knows that  

\begin{eqnarray*}
\sigma \left( \widehat{I}_{m,r} \right) &=& \dfrac{E_{m,r}}{E_{0,r}}\\
 \sigma\left( D \right)&=& E^p_{0,R}(X) \prod_{r=1}^{R-1} E^2_{0,r}(X) 
\end{eqnarray*}
where $p = 1$ for $n$ even and $p = 2$ for $n$ odd. If \eqref{symbol} were an algebra homomorphism, then the symbol of the product would be the product of the symbols and in our case it is evident that this product of the symbols is a polynomial and therefore non-singular. However $\sigma$ is only a linear bijection and the symbol of the product only equals the product of the symbols modulo terms of lower order. However, restricted to {\em invariants} it is known that a modification of $\sigma$ does realize an algebra isomorphism \cite{bib:duflo}. Whether or not something like this could be extended to {\em semi-invariants} of the type we consider here is currently under consideration. That is something that could also help with the issues raised earlier in this Section. 

\section*{Acknowledgements}  I wish to take this opportunity to express my deep and abiding gratitude to Hermann Flaschka for introducing me, many years ago, to the fascinating connections between Lie theory and Poisson geometry, as revealed by the Toda lattice. Hermann himself viewed the Toda lattice as a ``vehicle'' for exploring many fascinating areas of mathematics, which often had unexpected and beautiful relations with this nonlinear system. The work presented in this paper, motivated by our earlier study \cite{bib:efs} with Stephanie Singer, had its inception in extended explorations that Hermann and I undertook starting from early 2019 up to the time of his sudden passing in early 2021. I hope he would have liked this article, and that the ideas discussed here might help continue the mathematical legacy Hermann set in motion with his early seminal work on the Toda lattice. 

\appendix

\section{Appendix: Enveloping Algebras} \label{app}

We recall some basic definitions and facts. Associated to any Lie algebra, $(\frak{g}, [ \cdot , \cdot ])$ one has its {\em universal enveloping algebra} $\frak{U}(\frak{g})$ defined as the associative algebra 
\[
\frak{U}(\frak{g}) = T(\frak{g})/I
\]
where $T(\frak{g})$ is the tensor algebra of $\frak{g}$ and $I$ is the two-sided ideal generated by $\left\{ X \otimes Y - Y \otimes X - [X,Y] \right\}$ for $X,Y \in \frak{g}$. 
$\frak{g}$ is naturally included in $\frak{U}(\frak{g}) $ through $T^1(\frak{g})$. This has an analytical interpretation in terms of left-invariant differential operators. Given $X \in \frak{g}$, one defines a first order operator acting on smooth functions, $f$ on the associated group manifold $G$ by
\[
\tilde{X}f(g) = \frac{d}{dt} f(g \exp tX)|_{t=0}.
\]
(In our case we only be considered exponential groups so that $G = \exp(\frak{g}$. ) For any diffeomorphism $\gamma: G \to G$, one defines an action of $\gamma$ on $\tilde{X}$ by $\tilde{X}^\gamma f(g) = X f(g \circ \gamma) \circ \gamma^{-1} $. $\tilde{X}$ is left-invariant in the sense that taking $\gamma_h$ to be the diffeomorphism of left translation on $G$ by a fixed element $h \in G$ one has $\tilde{X}^{\gamma_h} = \tilde{X}$ for all such $h$. 

This extends naturally to monomial elements of $\frak{U}(\frak{g})$ by 
\[
(X_1 X_2 \dots X_k)^{\sim} \cdot f = \frac{\partial^k}{\partial t_1 \partial t_2 \cdots \partial t_k} f(g \exp t X_1 \exp t X_2 \dots \exp t X_k)_{t_1 = \cdots = t_k =0}
\]
and then to all of $\frak{U}(\frak{g}) $ by linearity. In fact one has an isomorphism of associative algebras \cite{bib:helg1},
\[
\frak{U}(\frak{g}) \simeq \mathbb{D}(G),
\]
where $\mathbb{D}(G)$ denotes the space of all left-invariant differential operators, $D$ on $G$ (meaning $D^{\gamma_h} = D$). Along these same lines the adjoint action of $G$ on  $\frak{g}$ extends to $\frak{U}(\frak{g})$,
\[
Ad_g(X_1 X_2 \dots X_k) = (Ad_g X_1) (Ad_g X_2) \dots (Ad_g X_k) 
\]
which corresponds to an adjoint action, $Ad_{g^{-1}}D$ on left-invariant differential operators. (Since $D$ is left-invariant this amounts to the action induced through right translation by $g$.) One then defines
\[
ad_X D \doteq \frac{d}{dt} Ad_{\exp(tX)} D|_{t=0} 
\]
which coincides with $ad$ on $\frak{g}$. It then follows from the product rule that
\[
ad_Y  X_1 \dots X_k = Y X_1 X_2\dots X_k  - X_1 X_2\dots X_k Y.
\]
It follows that if an element of $\frak{U}(\frak{g}) $ commutes with all elements of $\frak{g}$ then it is an element of the center of $\frak{U}(\frak{g}) $ denoted 
$\frak{Z}(\frak{g}) $ which in turn is isomorphic to the commutative algebra of bi-invariant differential operators.

Having made these identifications one is led to %an alternative representation that provides a more 
a direct connection to invariant theory. 
%Let 
%$\{X_1, \dots X_n \}$ be a basis for $\frak{g}$ which one may specifically regard as a basis for $T_e G$. For the exponential groups we consider,
%$(t_1, \dots, t_n) \mapsto \exp(t_1 X_1 + \cdots + t_n X_n)$ provides a global coordinate system on $G$. 
Let $\mathbb{S}(\frak{g})$ be the symmetric algebra over 
$\frak{g}$. Then \cite{bib:helg2} $\exists !$ linear bijection 
\begin{eqnarray} \label{symbol}
\sigma : \mathbb{S}(\frak{g}) &\to& \mathbb{D}(G)\\ \nonumber
\ni \sigma(X^m) &=& \tilde{X}^m.
\end{eqnarray}
This is not an algebra isomorphism; however, $\sigma$ commutes with the $Ad_G$ action on $\mathbb{S}(\frak{g})$ (which preserves $\mathbb{S}(\frak{g})$). It is natural to identify 
$\mathbb{S}(\frak{g})$ with the space of polynomial functions on $\frak{g}^*$ Thus $Ad_G$-invariant
polynomials on $\frak{g}$ map to $\frak{Z}(\frak{g}) $. Moreover, if the polynomials $\{ P_1, \dots, P_m\}$ generate the subalgebra of invariants, $\mathbb{S}(\frak{g})^G$, then $\{ \sigma(P_1), \dots, \sigma(P_m) \}$ generate $\frak{Z}(\frak{g}) $ as an algebra of bi-invariant operators. 

%% ================================================================

%% ================================================================


\begin{thebibliography}{}
%\bibitem[Ad79]{bib:adler}
%M. Adler.
%\newblock {\em On a trace functional for formal pseudodifferential operators and the symplectic structure for the Korteweg-de Vries type equations}.
%\newblock \textit{Invent. Math.} \textbf{50}, 219-248 (1979).

%\bibitem[ADH17]{bib:adh}
%A. Auffinger, M. Damron, and J. %Hanson.
%\newblock {\em 50 years of %first-passage percolation}.
%\newblock University Lecture Series, %vol 68. American Mathematical %Society, Providence, RI, 2017.


%\bibitem[Ai07]{bib:aigner}
%M. Aigner.
%\newblock {\em A Course in Enumeration}.
%\newblock Graduate texts in mathematics, vol 238. Springer, Berlin.

\bibitem[Arh79]{bib:arh}
A.A. Arhangel'skii. 
\newblock {\em Completely Integrable Hamiltonian Systems on the Group of Triangular Matrices}.
\newblock (Russian) Mat. Sb. (N.S.) \textbf{108 (150)} 134-142  (1979). 


%\bibitem[Ar07]{bib:arnold}
%V. I. Arnold.
%\newblock {\em Mathematical Methods of Classical Mechanics}.
%\newblock Graduate texts in mathematics, vol 60. Springer, Berlin.

\bibitem[AK71]{bib:ak}
L. Auslander and B. Kostant.
\newblock {\em Polarization and Unitary Representations of Solvable Lie Groups}.
\newblock \textit{Invent. Math.} \textbf{14}, 255-354 (1971).

%\bibitem[BFZ96]{bib:bfz}
%A. Berenstein, S. Fomin, and A. Zelevinsky.
%\newblock {\em Parametrizations. of Canonical Bases and Totally Positive Matrices}.
%\newblock \textit{Advances in  Mathematics} \textbf{122}, 49-149 (1996). 

\bibitem[BBO05]{bib:bbo}
P. Biane, P. Bougerol, and N. O'Connell.
\newblock {\em Littleman Paths and Brownian Paths}.
\newblock \textit{Duke Math. J.} \textbf{130}, 127-167, (2005). 

\bibitem[BGR17]{bib:bgr}
A.M. Bloch, F. Gay-Balmaz, and T.S. Ratiu.
\newblock {\em The Geometric Nature of the Flaschka Transformation}.
\newblock \textit{Communications in Mathematical Physics} \textbf{352}, 457-517 (2017). 

\bibitem[BBF11]{bib:bbf}
B. Brubaker, D. Bump and  S. Friedberg.
\newblock {\em Eisenstein Series, Crystals, and Ice}.
\newblock \textit{Notices Amer. Math. Soc.} \textbf{58}, 1563-1571, (2011). 


%\bibitem[CS21]{bib:cs}
%D. Croydon and M. Sasada.
%\newblock {\em Generalized hydrodynamic limit for the box-ball system}.
%\newblock Communications in Mathematical Physics, \textbf{383}, 427-463 (2021).


%\bibitem[DLT89]{bib:dlt}
%P. Deift, L. C. Li, and C. Tomei.
%\newblock {\em Matrix Factorizations and Integrable Systems}.
%\newblock \textit{Communications on Pure and Applied Mathematics} \textbf{42}, 443-521 (1989). 

%\bibitem[DNT83]{bib:dnt}
%P.Deift, T. Nanda, and C. Tomei.
%\newblock {\em Ordinary differential equations and the symmetric eigenvalue problem}.
%\newblock SlAM {\em Jl numer. Analysis} \textbf{20}, 1-22 (1983). 

\bibitem[DLNT86]{bib:dlnt}
P.Deift, L.C. Li, T. Nanda, and C. Tomei.
\newblock {\em The Toda flow on a generic orbit is integrable}.
\newblock \textit{Communications on Pure and Applied Mathematics} \textbf{39}, 183-232 (1986).

\bibitem[Dix59]{bib:dix}
J. Dixmier.
\newblock {\em Sur les Representations Unitaires des Groupes de Lie Nilpotents IV.}
\newblock Canadian Joournal of Mathematics, \textbf{11}, 321-344  (1959).

\bibitem[Du77]{bib:duflo}
M. Duflo.
\newblock {\em Op\'erateurs Diff\'erentiels Bi-invariants sur un Groupe de Lie.}
\newblock Ann. Sci. \'Ecole Norm. Sup., \textbf{10}, 265-288  (1977).

\bibitem[DM72]{bib:dm}
H. Dym and H. P. McKean.
\newblock {\em Fourier Series and Integrals}.
\newblock Probability and Mathematical Statistics,  {\bf 14},  Academic Press,  New York  (1972).

%\bibitem[EFH91]{bib:efh}
%N. M. Ercolani, H. Flaschka, and L. Haine.
%\newblock {\em Painlev\'{e} Balances and Dressing Transformations}.
%\newblock In: Painlev\'{e} Transcendents, NATO ASI series, Series B, Physics 278 (1991).

\bibitem[EFS93]{bib:efs}
N. M. Ercolani, H. Flaschka, and S. Singer.
\newblock {\em The Geometry of the Full Kostant-Toda Lattice}.
\newblock In {\em{Integrable systems (Luminy, 1991)}}, vol. 115 of {\em{Progr. Math.}} Birkh\"{a}user Boston, Boston, MA. 181-225 (1993).

%\bibitem[EM01]{bib:em}
%N. M. Ercolani and K. D. T.-R. McLaughlin, {\em Asymptotics and Integrable Structures for Biorthogonal Polynomials Associated to a Random Two-Matrix Model},  %\newblock {Physica D}, \textbf{152-153}, 232-268 (2001).

%\bibitem[ER21a]{bib:era}
%N. M. Ercolani and J. Ramalheira-Tsu.
%\newblock {\em The Ghost-Box-Ball System: A Unified Perspective on Soliton Cellular Automata, the RSK Algorithm and Phase Shifts}.
%\newblock {Physica D}, \textbf{426} (2021) 132986

%\bibitem[ER21b]{bib:er}
%N. M. Ercolani and J. Ramalheira-Tsu.
%\newblock {\em A Path-Counting Analysis of Phase Shifts in Box-Ball Systems}.
%\newblock 	arXiv:2106.07129 (2021).

%\bibitem[EW19]{bib:ew}
%N. M. Ercolani and P. Waters.
%\newblock {\em Relating Random Matrix Map enumeration to a universal symbol calculus for recurrence operators in terms of Bessel-Appel polynomials}.
%\newblock arXiv:1907.08026 (to appear in Random Matrices: Theory and Applications)

%\bibitem[FIYIN13]{bib:fiyin}
%A. Fukuda, E. Ishiwata, Y. Yamamoto, M. Iwasaki, and Y. Nakamura.
%\newblock {\em Integrable discrete hungry systems and their related matrix eigenvalues}.
%\newblock {\em{Annal. Mat. Pura Appl.}}, Vol \textbf{192}, 423-445 (2013).

\bibitem[Fl74]{bib:fl}
H. Flaschka.
\newblock {\em The Toda lattice II. Existence of Integrals}.
\newblock Physical Review B, \textbf{9}, 1924-1925  (1974).

%\bibitem[FH91]{bib:fh}
%H. Flaschka and L. Haine.
%\newblock {\em Vari\'{e}t\'{e}s de drapeaux et r\'{e}seaux de Toda}.
%\newblock Math. Z. \textbf{208}, 545-556 (1991).

\bibitem[Fo89]{bib:foll}
G.B. Folland.
\newblock {\em Harmonic Analysis in Phase Space}.
\newblock ,  {\bf 129},  Princeton University Press  (1989).

\bibitem[FuHa91]{bib:fuha}
W. Fulton and J. Harris.
\newblock {\em Representation Theory: A First Course}.
\newblock Annals of Mathematics Studies,  {\bf 122},  Springer-Verlag,  New York  (1991).

%\bibitem[FZ00]{bib:fz}
%S. Fomin and A. Zelevinsky.
%\newblock {\em Total Positivity: Tests and Parametrizations}.
%\newblock Math Intelligencer \textbf{22}, 23-33 (2000).

%\bibitem[Ga22]{bib:ganguly}
%S. Ganguly.
%\newblock {\em Random Metric Geometries on the Plane and Karder-Parisi-Zhang Universality}.
%\newblock Notices of the AMS \textbf{69}(1), 26 - 35 (2022).

%\bibitem[GGMS87]{bib:ggms}
%I.M. Gelfand, R.M. Goresky, R.D. MacPherson and V.V. Serganova.
%\newblock {\em Combinatorial Geometries, Convex Polyhedra, and Schubert Cells}.
%\newblock {Advances in Mathematicds}, \textbf{63}, 301-316 (1987). 

\bibitem[GS99]{bib:gs}
M. I. Gekhtman and M. Z. Shapiro.
\newblock {\em Non-Commutative and Commutative Integrability of Generic Toda Flows in Simple Lie Algebras}.
\newblock \textit{Communications on Pure and Applied Mathematics}, \textbf{52}(1), 53 - 84 (1999). 

\bibitem[Gi97]{bib:gi}
A. Givental.
\newblock {\em Stationary Phase Integrals, Quantum Toda Lattices, Flag Manifolds and the Mirror Conjecture}.
\newblock In {\em Topics in Singularity Theory}, American Mathematical Society Translations, Ser. 2, AMS, Providence, RI (1997).

\bibitem[GuSt84]{bib:gust}
V. Guillemin and S. Sternberg.
\newblock {\em Symplectic Techniques in Physics}.
\newblock Cambridge University Press,  (1984).

\bibitem[He78]{bib:helg1}
S. Helgason.
\newblock {\em Differential Geometry, Lie Grooups and Symmetric Spaces}.
\newblock AMS,  (1978).

\bibitem[He84]{bib:helg2}
S. Helgason.
\newblock {\em Grooups and Geometric Analysis}.
\newblock AMS,  (1984).

%\bibitem[Hi77]{bib:h}
%R. Hirota.
%\newblock {\em Nonlinear Partial Difference Equations. II. Discrete-Time Toda Equation}.
%\newblock J. Phys. Soc. Japan \textbf{43}(6), 2074-2078 (1977).

%\bibitem[Hi81]{bib:hirota81} 
%\newblock {\em Discrete analogue of a generalized Toda equation}.
%\newblock J. Phys. Soc. Japan., \textbf{50}, 3785-3791 (1981).

%\bibitem[HT95]{bib:ht}
%R. Hirota and S. Tsujimoto
%\newblock {\em Conserved Quantities of a Class of Nonlinear Difference-Difference Equations}.
%\newblock J. Phys. Soc. Japan \textbf{64}, 3125-3127 (1995).

%\bibitem[Ho92]{bib:howe}
%R. Howe.
%\newblock {\em A Century of Lie Theory}.
%\newblock In: Mathematics into the Twenty-first Century, AMS, (1992).

\bibitem[Jo77]{bib:jo}
A.A. Joseph.
\newblock {\em A PreparationTheorem for the Prime Spectrum of a Semisimple Lie Algebra}.
\newblock Journal of Algebra \textbf{48}, 241-289 (1977).


%\bibitem[KNW09]{bib:knw}
%S. Kakei, J. Nimmo, and R. Willox.
%\newblock {\em Yang–Baxter maps and the discrete KP hierarchy}.
%\newblock { Glasg. Math. J.}, \textbf{51}, 107-119 (2009).

\bibitem[Ki62]{bib:kir}
A.A. Kirillov.
\newblock {\em Unitary Representations of nilpotent Lie Groups}. 
\newblock {Russian Math Surveys} \textbf{17}, 57 - 110 (1962).

\bibitem[Ki99]{bib:kir2}
A.A. Kirillov.
\newblock {\em Merits and Demerits of the Orbit Method}. 
\newblock {Bulletin of the AMS} \textbf{36}, 433 - 488 (1999).


%\bibitem[KW15]{bib:kw}
%Y. Kodama and L. Williams.
%\newblock {\em The Full Kostant–Toda Hierarchy on the Positive Flag Variety}.
%\newblock {Communications in Mathematical Physics}, \textbf{335}, 247-€"283 (2015).  

%\bibitem[KS18]{bib:ks}
%Y. Kodama and B. Shipman.
%\newblock {\em Fifty years of the finite nonperiodic Toda lattice: a geometric and topological viewpoint}.
%\newblock {Journal of Physics A: Mathematical and Theoretical}, \textbf{51} (2018) 353001

\bibitem[Ko70]{bib:kostant3}
B. Kostant.
\newblock {\em Quantization and Unitary Representations}. 
\newblock {Lecture Notes in Math.} \textbf{170}, 87 - 208 (1970).

\bibitem[Ko78]{bib:kostant}
B. Kostant.
\newblock {\em On Whittaker Vectors and Representation Theory}. 
\newblock {\em{Inventiones Math.}} \textbf{48}, 101-184 (1978).

\bibitem[Ko79]{bib:kostant2}
B. Kostant.
\newblock {\em The Solution to a Generalized Toda Lattice and Representation Theory}. 
\newblock {\em{Advances in Mathematics}} \textbf{34}, 195-338 (1979).

\bibitem[Ko12]{bib:kostant4}
B. Kostant.
\newblock {\em The Cascade of Orthogonal Roots and the Coadjoint Structure of the Nilradical of a Borel Subgroup of a Semisimple Lie Group}. 
\newblock {\em{Moscow Math. Journal}} \textbf{12}, 1-16 (2012).

\bibitem[Ko13]{bib:kostant5}
B. Kostant.
\newblock {\em Center $U(\frak{n})$, Cascade of Orthogonal Roots, and a Construction of Lipsman-Wolf}.
\newblock In {\em{Li Groups: Structure, Actions and Representations, In Honor of Joseph A. Wolf on the Occasion of his 75th Birthday}}, vol. 306 of 
{\em{Progr. Math.}} Springer Science+Business Media, New York, NY 163 - 173 (2013).

\bibitem[Kn86]{bib:knapp}
A.Knapp
\newblock {\em Representation Theory of Semisimple Groups}.
\newblock  Princeton University Press (1986).

\bibitem[LW78]{bib:lw}
R.L. Lipsman and J.A. Wolf.
\newblock {\em The Plancherel Formula for Parabolic Subgroups of the Classical Groups}. 
\newblock {Journal D'Analyse Mathe\'matique} \textbf{34}, 120 - 161 (1978).


%\bibitem[Li07]{bib:l}
%G. L. Litvinov.
%\newblock {\em The Maslov dequantization, idempotent and tropical mathematics: A brief introduction}.
%\newblock Journal of Mathematical Sciences, \textbf{140}(3), 426-444 (2007).

%\bibitem[LMRS11]{bib:lmrs}
%G. L. Litvinov, V. P. Maslov, A. YA. Rodionov, and A. N. Sobolevski.
%\newblock {\em Universal Algorithms, Mathematics of Semirings and Parallel Computations}.
%\newblock Lect. Notes Comput. Sci. Eng., \textbf{75}, 63–89 (2011).

%\bibitem[LS19]{bib:scrimshaw}
%X. Liu and T. Scrimshaw.
%\newblock {\em A uniform approach to soliton cellular automata using rigged configurations}.
%\newblock Ann. Henri Poincar\'{e}, \textbf{20}(4), 1175–1215 (2019).

%\bibitem[Lu94]{bib:lu}
%G. Lusztig.
%\newblock {\em Total positivity in reductive groups}.
%\newblock In {\em{Lie Theory and Geometry: In Honor of Bertram Kostant}}, vol. 123 of {\em{Progr. Math.}} Birkh\"{a}user Boston, Boston, MA, (1994).

\bibitem[Ma58]{bib:mackey}
G. Mackey.
\newblock {\em Unitary Representations of Grooup Extensions I}.
\newblock Acta Math., \textbf{99}, 265 - 311 (1958).

%\bibitem[Mo75]{bib:moser}
%J. Moser.
%\newblock {\em Finitely many mass points on the line under the influence of an exponential potential -- an integrable system}.
%\newblock Dynamical Systems, Theory and Applications, Lecture Notes in Physics \textbf{38}, Springer (1975).

\bibitem[MW73]{bib:mw}
C.C. Moore and J.A. Wolf.
\newblock {\em Square Integrble Representations of nilpotent Groups}. 
\newblock {Transactions of the AMS} \textbf{185}, 445 - 462 (1973).

%\bibitem[NY04]{bib:ny}
%M. Noumi and Y. Yamada.
%\newblock {\em Tropical Robinson-Schensted-Knuth correspondence and birational Weyl group actions}.
%\newblock Representation Theory of Algebraic Groups and Quantum Groups, \textbf{40}, 371-442 (2004).

\bibitem[O13]{bib:o}
N. O'Connell.
\newblock {\em Geometric RSK and the Toda lattice}.
\newblock Illinois Journal of Mathematics , \textbf{57}(3), 883-918 (2013).

\bibitem[Pu71]{bib:puk}
L. Pukanszky.
\newblock {\em On the Theory of Exponential Groups}.
\newblock Trans. Amer. Math. Soc., \textbf{126}, 487-507 (1967).

\bibitem[Ra20]{bib:r}
J. Ramalheira-Tsu.
\newblock {\em The Kostant-Toda Lattice, Combinatorial Algorithms and Ultradiscrete Dynamics}.
\newblock Available from ProQuest Dissertations \& Theses Global (2461615150). (2020).

\bibitem[Sh10]{bib:sh}
F. Shahidi.
\newblock {\em Eisenstein Series and Automorphic L-functions}.
\newblock American Mathematical Society Colloquium Publications,  Providence, RI  (2010).

%\bibitem[SIK20]{bib:sik}
%M. Shinjo, M. Iwasaki and K. Kondo.
%\newblock {\em The Kostant-Toda equation and the hungry integrable systems}.
%\newblock  Journal of Mathematical Analysis and Applications, \textbf{483},  (2020) 123627

\bibitem[Si05]{bib:singer}
S.F. Singer.
\newblock {\em Linearity, Symmetry and Prediction in the Hydrogen Atom}.
\newblock Undergraduate Texts in Mathematics,  Springer-Verlag,  New York  (2005).


%\bibitem[Sp21]{bib:spohn}
%H. Spohn.
%\newblock {\em Hydrodynamic Equations for the Toda Lattice}.
%\newblock arXiv:2101.06528

%\bibitem[St93]{bib:strang}
%G. Strang.
%\newblock {\em Introduction to linear algebra}.
%\newblock Volume 3. Wellesley-Cambridge Press Wellesley, MA, (1993).

%\bibitem[Su18]{bib:suris}
%Y. B. Suris.
%\newblock {\em Discrete time Toda systems}.
%es78\newblock J. Phys. A: Math. Theor. \textbf{51}(33), Special Issue ``Fifty Years of the Toda lattice'' (2018).

%\bibitem[Sy78]{bib:symes78}
%W. W. Symes.
%\newblock {\em Systems of Toda type, inverse spectral problems and representation theory.}
%\newblock \textit{Invent. Math.}, \textbf{59}, 13-53 (1978).

\bibitem[Sy80]{bib:symes}
W. W. Symes.
\newblock {\em Hamiltonian group actions and integrable systems}.
\newblock \textit{Physica D}, \textbf{1}, 339-374 (1980).

\bibitem[Th82]{bib:thimm}
A. Thimm.
\newblock {\em Integrable geodesic flows on homogeneous spaces}.
\newblock Ergod. Th. \& Dynam. Sys. \textbf{1}, 495-517 (1982).

\bibitem[To67]{bib:toda}
M. Toda.
\newblock {\em Vibration of a chain with a non-linear interaction}.
\newblock J. Phys. Soc. Jpn., \textbf{22}(2), 431-436 (1967).

%\bibitem[To04]{bib:tokihiro}
%T. Tokihiro.
%\newblock {\em Ultradiscrete Systems (Cellular Automata)}.
%\newblock Lect. Notes in Phys., \textbf{644}, 383-424 (2004).

%\bibitem[TS90]{bib:ts}
%D. Takahashi, and J. Satsuma.
%\newblock {\em A Soliton Cellular Automaton}.
%\newblock J. Phys. Soc. Jpn., \textbf{59}, 3514-3519 (1990).

%\bibitem[TTS96]{bib:tts}
%M. Torii, D. Takahashi, and J. Satsuma.
%\newblock {\em Combinatorial representation of invariants of a soliton cellular automaton}.
%\newblock Physica D, \textbf{92}, 209-220 (1996).

%\bibitem[TNS99]{bib:tns}
%T. Tokihiro, A. Nagai, and J. Satsuma.
%\newblock {\em A Proof of solitonical nature of box and ball systems by means of inverse ultra-discretization}.
%\newblock Inverse Problems, \textbf{15}, 1639-1662 (1999).

\bibitem[Tro80]{bib:trof}
V.V. Trofimov.
\newblock {\em Finite-Dimensional Representations of Lie Algebras and Completely Integrable Systems.} 
\newblock (Russian) Mat. Sb. (N.S.) \textbf{111 (153)} 610-621  (1980).

%\bibitem[Vi01]{bib:v}
%O. Viro.
%\newblock {\em Dequantization of real algebraic geometry on a logarithmic paper}.
%\newblock in: 3rd European Congress of Mathematics, Barcelona 2000, Vol. I, Birkh\"{a}user, Basel, 135-146 (2001).

%\bibitem[Wa84]{bib:watkins}
%D. Watkins.
%\newblock {\em Isospectral Flows}.
%\newblock SIAM Review, \textbf{26}, 379-391 (1984).

\bibitem[Wo14]{bib:wo}
J.A. Wolf.
\newblock {\em Plancherel Formula for Minimal Parabolic Subgroups}. 
\newblock {\em{Journal of Lie Theory}} \textbf{24}, 791 - 808 (2014).


\bibitem[Wo16]{bib:wo1}
J.A. Wolf.
\newblock {\em Stepwise Square Integrable Representations: The Concept and Some Consequences}. 
\newblock  in: Lie Theory and its Applications in Physics, Varna, Bulgaria, June 2015, 
Springer Proceedings in Mathematics \& Statisitics, Vol. 191, Springer Nature, Singapore, 181-202 (2016).


%\bibitem[YF10]{bib:yf}
%Y. Yamamoto and T. Fukaya.
%\newblock {\em Differential qd algorithm for totally nonnegative Hessenberg matrices: introduction of origin shifts and relationship with the discrete hungry Lotka-%Volterra system}.
%\newblock {\em JSIAM Lett.}, \textbf{2}, 69-72 (2010).

\end{thebibliography}
\end{document}